\documentclass{amsart}
\usepackage[utf8]{inputenc}
\usepackage{amsmath}
\usepackage{amssymb}
\usepackage{leftidx}
\usepackage{amsthm}
\usepackage{hyperref}
\usepackage{tikz}
\usetikzlibrary{intersections}
\usepackage{standalone}
\usepackage{tkz-euclide}
\usepackage{enumerate}

\newcommand{\R}{\mathbb{R}}
\newcommand{\N}{\mathbb{N}}

\newcommand{\Z}{\mathbb{Z}}
\newcommand{\T}{\mathbb{T}}
\newcommand{\C}{\mathbb{C}}
\newcommand{\ofrac}[2]{\frac{#1}{#2}\oslash}
\newcommand{\dcm}{\Omega}
\newcommand{\dcmg}{h}
\newcommand{\dcn}{J}

\newtheorem{aptheorem}{Theorem}[section]
\newtheorem{aplemma}[aptheorem]{Lemma}

\newtheorem{apcorollary}[aptheorem]{Corollary}
\newtheorem{apconjecture}[aptheorem]{Conjecture}

\theoremstyle{definition}

\newtheorem{definition}[aptheorem]{Definition}
\newtheorem{example}[aptheorem]{Example}

\newtheorem{apremark}[aptheorem]{Remark}

\numberwithin{equation}{section}

\begin{document}


\title{Tropical second main theorem and the Nevanlinna inverse problem}
\author{Juho Halonen, Risto Korhonen and Galina Filipuk}
\date{\today}
\thanks{GF acknowledges the support of the grant entitled ”Geometric approach to ordinary differential equations” funded under New Ideas 3B competition within Priority Research Area III implemented under the “Excellence Initiative – Research University” (IDUB) Programme (University of Warsaw) (nr 01/IDUB/2019/94). GF is also partially supported by the Ministry of Science and Innovation of Spain and the European Regional Development Fund (ERDF) grant number PID2021-124472NB-I00.}

\address{Juho Halonen, Department of Physics and Mathematics, University of Eastern Finland, P.O. Box 111 FI-80101 Joensuu, Finland}
\address{Risto Korhonen, Department of Physics and Mathematics, University of Eastern Finland, P.O. Box 111 FI-80101 Joensuu, Finland}
\address{Galina Filipuk, Institute of Mathematics, University of Warsaw, ul. Banacha 2, 02-097, Warsaw, Poland}
\address{email addresses: juho.halonen@uef.fi, risto.korhonen@uef.fi, filipuk@mimuw.edu.pl}


\maketitle

\begin{abstract}
A generalization of the second main theorem of tropical Nevanlinna theory is presented for noncontinuous piecewise linear functions and for tropical hypersurfaces without requiring a growth condition. The method of proof is novel and significantly more straightforward than previously known proofs. The tropical analogue of the Nevanlinna inverse problem is formulated and solved for tropical meromorphic functions and tropical hypersurfaces.

\end{abstract}

\section{Introduction}

Nevanlinna's second main theorem is a deep generalization of Picard's theorem, providing a quantitative estimate of how frequently a meromorphic function takes on a given value, in terms of the growth of the function \cite{nevanlinna:25}. This quantification can be expressed explicitly as the deficiency relation
$$
\sum_a \delta(a,f) \leq 2 ,
$$
where the sum is taken over all $a$ in the Riemann sphere and $\delta(a,f)$ is the Nevanlinna defect of a meromorphic function $f$ at the point $a$. The Nevanlinna inverse problem, on the other hand, involves the construction of meromorphic functions that have pre-assigned deficiencies at specified points (see the appendix for more details on the inverse problem). This problem was solved by Drasin using quasiconformal mappings \cite{drasin:77}. In this paper we prove, on one hand, a tropical analogue of the second main theorem for piecewise linear functions without growth restrictions, and, on the other hand, give a complete solution to the tropical Nevanlinna inverse problem.  

The tropical semiring is the set of real numbers extended with an element $-\infty$, denoted as $\T = \R\cup\{-\infty\}$, equipped with two operations: $\oplus$, the maximum operation defined as $a\oplus b = \max\{a,b\}$, and $\otimes$, the usual addition of real numbers, defined as $a\otimes b =  a + b$ \cite{maclagans:15}.
The additive and multiplicative inverses of the tropical semiring are $0_\circ = -\infty$ and $1_\circ = 0$. The tropical semiring is categorized as a semiring since not all of its elements possess an additive inverse. A solution cannot be found for equations such as $5\oplus x = 0_\circ$, which means that subtraction cannot be performed within the tropical semiring. Division in the tropical semiring is defined as $a\oslash b = a - b$ and exponentiation as $a^{\otimes\alpha} = \alpha a$. 

    A continuous piecewise linear function $f:\R\to\R$ is said to be tropical meromorphic.
A point of discontinuity $x$ of the derivative of a tropical meromorphic function $f$ is called a pole, when
\begin{equation*}
    \omega_f(x) := \lim_{\varepsilon\to0^+}(f'(x+\varepsilon)-f'(x-\varepsilon)) < 0
\end{equation*}
and a root when $\omega_f(x) > 0$. The multiplicity of a root or a pole at $x$ is $\tau_f(x) := |\omega_f(x)|$ \cite{halburds:09}.

Tropical Nevanlinna theory studies the growth and complexity of tropical meromorphic functions in an analogous way the classical Nevanlinna theory studies complex meromorphic functions. Halburd and Southall developed tropical Nevanlinna theory initially as a means to examine the integrability of ultra-discrete equations~\cite{halburds:09}. By utilizing tropical Nevanlinna theory, they proposed that there is a connection between the existence of finite-order tropical meromorphic solutions for ultra-discrete equations and an ultra-discrete version of the Painlev\'e property. Halburd and Southall proved a number of central results of tropical Nevanlinna theory including a tropical analogue of the Poisson-Jensen formula, the first main theorem, a tropical analogue of the lemma on the logarithmic derivative for finite order tropical meromorphic functions and a tropical analogue of Clunie's lemma. Building on this foundation, Laine and Yang \cite{lainey:10} proceeded to prove the Mohon'ko and Valiron-Mohon'ko lemmas and also a stronger version of Clunie's lemma. Until that point, tropical Nevanlinna theory was limited to tropical meromorphic functions with only integer slopes. Laine and Tohge \cite{lainet:11} extended tropical Nevanlinna theory for arbitrary real valued slopes and noted that the results listed above still hold in the extended setting. They also proved a tropical second main theorem for tropical meromorphic functions with hyper-order less than $1$.

Korhonen and Tohge \cite{korhonent:16AM} extended tropical Nevanlinna theory for $n$-dimensional tropical projective space and proved a tropical analogue of Cartan's second main theorem, which implies the second main theorem by Laine and Tohge. Cao and Zheng \cite{caozheng:18} further extended tropical Nevanlinna theory for $n$-dimensional tropical hypersurfaces, and weakened the growth condition of the tropical analogue of the lemma on the logarithmic derivative from functions with hyper-order less than $1$ to subnormal functions. Their improvement of the growth condition in the tropical analogue of the lemma on the logarithmic derivative directly implies a corresponding  improvement in the growth conditions of the tropical second main theorem and many other results that rely on the lemma on the logarithmic derivative. For more background to the topics discussed in this paper, see the appendix at the end of the paper, and for a general overview on tropical Nevanlinna theory see \cite{korhonenlainetohge:15}. For more details about tropical mathematics in general see \cite{maclagans:15}.

Laine and Tohge \cite{lainet:11}  posited that the tropical second main theorem would not hold if the growth condition was omitted. In this paper, we will demonstrate that the growth condition can be completely excluded.

Halburd and Southall \cite{halburds:09} have shown that a large class of ultra-discrete equations, containing equations of Painlev\'e type, admits infinitely many continuous piecewise linear solutions. With the same method, but without the restriction of requiring continuity, it can be shown that the same class of equations also admits infinitely many noncontinuous piecewise linear solutions. In the first part of this paper we will extend tropical Nevanlinna theory for noncontinuous piecewise linear functions proving the Poisson-Jensen formula and the first and the second main theorems for these functions. The first and second main theorems will also be extended for a class of piecewise linear target values and the growth condition will be dropped from the second main theorem. The proof of the second main theorem will be greatly simplified compared to the previous proofs.

In the second part of the paper, we will prove an improved version of the second main theorem for tropical hypersurfaces, in which the growth condition will be dropped and the inequalities will be made tighter.

In the last part of this paper, we will formulate and answer three different versions of the inverse problem in the context of tropical Nevanlinna theory. First, we will show that for any $\delta\in[0,1]$, we can find a tropical rational function $f$ and a constant $a\in\R$ such that $\delta(a,f) = \delta$. Second, we will demonstrate that there exists a tropical meromorphic function $f$ such that for all $\delta\in[0,1]$, there exists $a\in\R$ such that $\delta(a,f) = \delta$. Lastly we will formulate and answer the inverse problem for tropical hypersurfaces. We will also examine some properties of the defect as a real valued function.  Finally we will disprove the tropical version of Griffiths conjecture \cite{griffiths:72} which was proposed by Cao and Zheng \cite[Conjecture 4.11]{caozheng:18}.

\section{The second main theorem for piecewise linear functions}

\begin{definition}
    Let $f:\R\to\R$. If there exists disjoint intervals $I_k\subset\R,\,k\in\N$ such that each of them contains more than one element,
    \begin{equation*}
        \bigcup_{k\in\N}I_k = \R
    \end{equation*}
    and for each interval $I_k$ there exists $\alpha,\beta\in\R$ such that $f(x) = \alpha x +\beta$ for all $x\in I_k$, then $f$ is said to be piecewise linear. 
\end{definition}

Halburd and Southall \cite{halburds:09} described the following method for generating continuous piecewise linear solutions for ultra-discrete equations of the type
\begin{equation}\label{equation: halburd ultra-discrete}
    y(x+1)\otimes y(x-1) = R(x,y(x)),
\end{equation}
where $R$ is a tropical rational function of $x$ and $y$. Choose any values for $y(0)$ and $y(1)$. Then compute $y(2) = R(1,y(1))-y(0)$. Now if you define $y(x)$ on the intervals $(0,1)\cup (1,2)$ in a way that $y(x)$ is a continuous piecewise linear function on the interval $[0,2]$, then the equation \eqref{equation: halburd ultra-discrete} extends $y$ uniquely to a tropical meromorphic function on the whole real line. However, if we allow discontinuities for the piecewise linear function in the initial interval $[0,2]$, then we can generate infinitely many noncontinuous piecewise linear solutions to the equation \eqref{equation: halburd ultra-discrete}. Motivated by this we will extend tropical Nevanlinna theory for noncontinuous piecewise linear functions.

Discontinuities are classified into three different categories: removable discontinuities, jump discontinuities and essential discontinuities \cite{apostol:74}. We say that a point of discontinuity $x_0$ of $f:\R\to\R$ is a jump discontinuity if the one-sided limits
\begin{equation*}
    \lim_{x\to x_0^+}f(x) =: f(x_0+)\in\R
\end{equation*}
and
\begin{equation*}
    \lim_{x\to x_0^-}f(x) =: f(x_0-)\in\R
\end{equation*}
exist and $f(x_0+)\neq f(x_0-)$. If the one-sided limits above exist, but $f(x_0+) = f(x_0-)$ then $x_0$ is said to be a removable discontinuity of $f$, and if one or both of the one-sided limits above do not exist in $\R$, then $x_0$ is said to be an essential discontinuity of $f$. 
\begin{aplemma}
    All discontinuities of a piecewise linear function are jump discontinuities.
\end{aplemma}
\begin{proof}
    Let $f:\R\to\R$ be a piecewise linear function. By definition there exists a partition of $\R$ into intervals $I_k$ such that $f$ is linear on each interval. If $f$ is discontinuous at $x_0$, then $x_0$ must be at the ends of two of the intervals. We know that $f$ is linear on both intervals and therefore the one-sided limits $f(x_0+)$ and $f(x_0-)$ must exist. Since $f$ is continuous on each of the intervals and $x_0$ is included in one of the intervals we must have either $f(x_0+) = f(x_0)$ or $f(x_0-) = f(x_0)$. On the other hand, since $f$ is discontinuous at $x_0$ we must also have either $f(x_0+) \neq f(x_0)$ or $f(x_0-) \neq f(x_0)$. Therefore we arrive at the conclusion that $f(x_0+) \neq f(x_0-)$.
\end{proof}
From the proof above we can also see that a piecewise linear function must be left-continuous or right-continuous at every point.

\begin{figure}[h!]
    \centering
    \begin{tabular}{c c c c}
        \begin{tikzpicture}
\def\xmin{-1}
\def\xmax{1}
\tkzInit[xmax=\xmax,ymax=1,xmin=\xmin,ymin=-1]
   \begin{scope}
    \end{scope}
   \tkzDrawX[noticks = true, label = ]
   \tkzDrawY[noticks = true, label = ]

\draw [domain=\xmin:0] plot(\x,1) node[right]{};
\draw [domain=0:\xmax] plot(\x,-1) node[right]{};
\filldraw[color = black, fill = white] (0,1) circle (2pt);
\filldraw (0,-1) circle (2pt);





\end{tikzpicture} & \begin{tikzpicture}
\def\xmin{-1}
\def\xmax{1}
\tkzInit[xmax=\xmax,ymax=1,xmin=\xmin,ymin=-1]
   \begin{scope}
    \end{scope}
   \tkzDrawX[noticks = true, label = ]
   \tkzDrawY[noticks = true, label = ]

\draw [domain=\xmin:0] plot(\x,-1) node[right]{};
\draw [domain=0:\xmax] plot(\x,1) node[right]{};
\filldraw[color = black, fill = white] (0,1) circle (2pt);
\filldraw (0,-1) circle (2pt);





\end{tikzpicture} & \begin{tikzpicture}
\def\xmin{-1}
\def\xmax{1}
\tkzInit[xmax=\xmax,ymax=1,xmin=\xmin,ymin=-1]
   \begin{scope}
    \end{scope}
   \tkzDrawX[noticks = true, label = ]
   \tkzDrawY[noticks = true, label = ]

\draw [domain=\xmin:0] plot(\x,1) node[right]{};
\draw [domain=0:\xmax] plot(\x,-1) node[right]{};
\filldraw[color = black, fill = white] (0,-1) circle (2pt);
\filldraw (0,1) circle (2pt);





\end{tikzpicture} & \begin{tikzpicture}
\def\xmin{-1}
\def\xmax{1}
\tkzInit[xmax=\xmax,ymax=1,xmin=\xmin,ymin=-1]
   \begin{scope}
    \end{scope}
   \tkzDrawX[noticks = true, label = ]
   \tkzDrawY[noticks = true, label = ]

\draw [domain=\xmin:0] plot(\x,-1) node[right]{};
\draw [domain=0:\xmax] plot(\x,1) node[right]{};
\filldraw[color = black, fill = white] (0,-1) circle (2pt);
\filldraw (0,1) circle (2pt);





\end{tikzpicture}\\
    \end{tabular}
    \caption{From left to right, two  negative jumps and two of positive jumps at $0$.}
    \label{fig:disc roots and poles at the origin}
\end{figure}
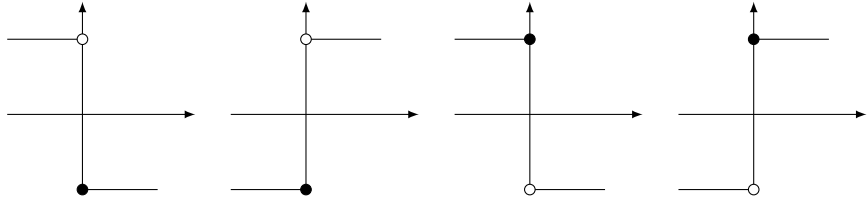

For a piecewise linear function $f$ we define
\begin{equation*}
    \dcm_f(x) := f(x+)-f(x-)
\end{equation*}
and we say that $x\neq 0$ is a positive jump if $x\dcm_f(x)>0$ and a negative jump if $x\dcm_f(x) < 0$. If there is a discontinuity at $0$, we say that it is a positive jump if $f(0) = \max\{f(0+),f(0-)\}$ and a negative jump if $f(0) = \min\{f(0+),f(0-)\}$ (see Figure \ref{fig:disc roots and poles at the origin}). Since all the discontinuities of a piecewise linear function are jump discontinuities, we know that $\dcm_f(x) = 0$ if and only if $f$ is continuous at $x$. 
The height of a jump at $x$ is $\dcmg_f(x) := |\dcm_f(x)|$. When moving away from the origin, a piecewise linear function jumps up at a positive jump and down at a negative jump. 

\begin{aplemma}\label{lemma: no accumulation points for discontinuities}
    Let $f$ be a piecewise linear function. Then the set $\{x\in\R:\dcm_f(x)\neq 0\}$ containing the points of discontinuity of $f$ cannot have accumulation points.
\end{aplemma}
\begin{proof}
    Let $f:\R\to\R$. Suppose that the set of discontinuities of $f$ has an accumulation point $x_0$. Then any interval of $\R$ that contains $x_0$ also contains another discontinuity. This means that on any interval that contains $x_0$,
    there do not exist $\alpha,\beta\in\R$ such that $f(x) = \alpha x + \beta$ on the whole interval. Therefore $f$ cannot be piecewise linear. 
\end{proof}
By Lemma \ref{lemma: no accumulation points for discontinuities} we can be sure, that on any finite interval a piecewise linear function has only finitely many discontinuities. 

The Poisson-Jensen formula for continuous piecewise linear functions was proven by Halburd and Southall.
\begin{aptheorem}[\cite{halburds:09}]\label{theorem: poisson jensen old}
    Suppose $f$ is a continuous piecewise linear function on $[-r,r]$ for some $r>0$ and denote the distinct roots, resp. poles of $f$ in this interval by $a_\mu$, resp. by $b_\nu$ with their corresponding multiplicities $\tau_f$. Then for any $x\in(-r,r)$ we have the tropical Poisson-Jensen formula
    \begin{align*}
    f(x)=&\frac{1}{2}(f(r)+f(-r)) + \frac{x}{2r}(f(r)+f(-r))\\
    & - \frac{1}{2r}\sum_{|a_\mu|<r}\tau_f(a_\mu)(r^2-|a_\mu-x|r-a_\mu x)\\
    &+\frac{1}{2r}\sum_{|b_\nu|<r}\tau_f(b_\nu)(r^2-|b_\nu-x|r-b_\nu x).
\end{align*}
In particular, the case $x = 0$ gives the tropical Jensen formula
\begin{equation*}
\begin{aligned}
    f(0) =\frac{1}{2}(f(r)+f(-r))
    - \frac{1}{2}\sum_{|a_\mu|<r}\tau_f(a_\mu)(r-|a_\mu|)+\frac{1}{2}\sum_{|b_\nu|<r}\tau_f(b_\nu)(r-|b_\nu|).
\end{aligned}
\end{equation*}
\end{aptheorem}

We will generalise the Poisson-Jensen formula for piecewise linear functions that may have discontinuities.

\begin{aptheorem}\label{theorem: poisson jensen extended}
    Suppose $f$ is a piecewise linear function on $[-r,r]$ for some $r>0$ and denote the distinct roots, resp. poles of $f$ in this interval by $a_\mu$, resp. by $b_\nu$ with their corresponding multiplicities $\tau_f$ and the distinct positive jumps, resp. negative jumps of $f$ in this interval by $\alpha_\mu$, resp. by $\beta_\nu$ with their corresponding heights $\dcmg_f$. Then for any $x\in(-r,r)$ we have the Poisson-Jensen formula
    \begin{equation}\label{equation: poisson jensen statement}
        \begin{aligned}
            f(x) 
            & = \frac{1}{2}(f(r)+f(-r)) + \frac{x}{2r}(f(r)-f(-r))\\
            & - \frac{1}{2r}\sum_{|a_\mu|<r}\tau_f(a_\mu)(r^2 - |a_\mu - x|r - a_\mu x)\\
            & + \frac{1}{2r}\sum_{|b_\nu|<r}\tau_f(b_\nu)(r^2 - |b_\nu - x|r - b_\nu x)\\
            & - \frac{x}{2r}\left(\sum_{-r\leq \beta_\nu\leq 0}\dcmg_f(\beta_\nu)+\sum_{0\leq \alpha_\mu\leq r}\dcmg_f(\alpha_\mu)\right)\\
            & + \frac{x}{2r}\left(\sum_{-r\leq \alpha_\mu\leq 0}\dcmg_f(\alpha_\mu)+\sum_{0\leq \beta_\nu\leq r}\dcmg_f(\beta_\nu)\right)\\
            & - \frac{1}{2}\left(\sum_{-r\leq \alpha_\mu\leq \min\{0,x\}}\!\!\!\!\!\!\!\!\!\!\!\!\dcmg_f(\alpha_\mu)+\sum_{\max\{x,0\}\leq \alpha_\mu\leq r}\!\!\!\!\!\!\!\!\!\!\!\!\dcmg_f(\alpha_\mu) - \sum_{ x\leq\alpha_\mu\leq 0}\!\!\!\dcmg_f(\alpha_\mu) - \sum_{0\leq\alpha_\mu\leq x}\!\!\!\dcmg_f(\alpha_\mu)\right) \\
            & + \frac{1}{2}\left(\sum_{-r\leq \beta_\nu\leq \min\{0,x\}}\!\!\!\!\!\!\!\!\!\!\!\!\dcmg_f(\beta_\nu)+\sum_{\max\{x,0\}\leq \beta_\nu\leq r}\!\!\!\!\!\!\!\!\!\!\!\!\dcmg_f(\beta_\nu) - \sum_{ x\leq\beta_\nu\leq 0}\!\!\!\dcmg_f(\beta_\nu) - \sum_{0\leq\beta_\nu\leq x}\!\!\!\dcmg_f(\beta_\nu)\right) \\
            & - \frac{x}{2r}\dcm_f(0) -\frac{1}{2}\left(A_f(x) + B_f(x)\right),
        \end{aligned}
    \end{equation}
    where
    \begin{equation*}
        A_f(x) = 
        \begin{cases}
            0, &\text{if }f \text{ is continuous at }x,\\
           \dcm_f(x), &\text{if }f \text{ is left-discontinuous at }x,\\
           -\dcm_f(x), &\text{if }f \text{ is right-discontinuous at }x
        \end{cases}
    \end{equation*}
    and
    \begin{equation*}
        B_f(x) = 
        \begin{cases}
            0, &\text{if }x = 0,\\
           -\dcm_f(0), &\text{if }x< 0,\\
           \dcm_f(0), &\text{if }x > 0.
        \end{cases}
    \end{equation*}
    If $x = 0$ we have the Jensen formula
    \begin{equation*}
        \begin{aligned}
            f(0)
            & = \frac{1}{2}(f(r)+f(-r))\\
             &- \frac{1}{2}\sum_{|a_\mu|<r}\tau_f(a_\mu)(r - |a_\mu|)
             + \frac{1}{2}\sum_{|b_\nu|<r}\tau_f(b_\nu)(r - |b_\nu|)\\
             &- \frac{1}{2}\sum_{|\alpha_\mu|\leq r}\dcmg_f(\alpha_\mu)
             +\frac{1}{2}\sum_{|\beta_\nu|\leq r}\dcmg_f(\beta_\nu).
        \end{aligned}
    \end{equation*}
\end{aptheorem}
\begin{example}\label{example: poisson jensen concrete}
    Let $r=4$ and define $f:[-4,4]\to\R$,
    \begin{equation*}
        f(x) = \begin{cases}
            -x-2, &-4\leq x < -2,\\
            \max\{-x-1,2x+2\}, &-2\leq x < 0,\\
            -x + 1, &0\leq x\leq 1,\\
            -\max\{x-4,2x- 6\}, &1 < x < 4,\\
            -1, &x=4.
        \end{cases}
    \end{equation*}
    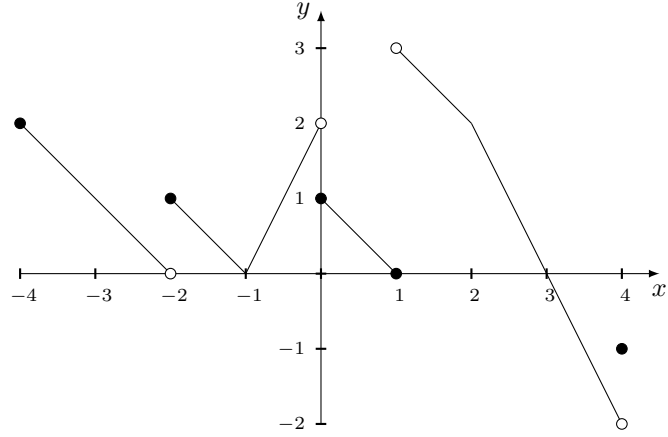
\begin{figure}
        \centering
        \begin{tikzpicture}
\def\xmin{-6}
\def\xmax{10}

\def\c{(-1,0)}
\def\ccv{(0,2)}
\def\cco{(0,1)}
\def\ccc{(2,2)}

\def\r{4}

\def\kv{(-2,0)}
\def\ko{(-2,1)}
\def\kkv{(1,0)}
\def\kko{(1,3)}
\def\kkkv{(4,-2)}
\def\kkko{(4,-1)}

\tkzInit[xmax=\r,ymax=3,xmin=-\r,ymin=-2]
   \begin{scope}

    \end{scope}

   \tkzLabelX[orig=false,step=1,node font=\scriptsize]
   \tkzLabelY[orig=false,step=1,node font=\scriptsize]
   \tkzDrawX
   \tkzDrawY

\draw (-4,2) -- \kv;
\draw \ko -- \c;
\draw \c -- \ccv;
\draw \cco -- \kkv;
\draw \kko -- \ccc;
\draw \ccc -- \kkkv;

\filldraw[color = black, fill = white] \kv circle (2pt);
\filldraw \ko circle (2pt);

\filldraw[color = black, fill = white] \ccv circle (2pt);
\filldraw \cco circle (2pt);

\filldraw[color = black, fill = white] \kko circle (2pt);
\filldraw \kkv circle (2pt);

\filldraw[color = black, fill = white] \kkkv circle (2pt);
\filldraw \kkko circle (2pt);

\filldraw (-4,2) circle (2pt);

\end{tikzpicture}
        \caption{Function $f$ in Example \ref{example: poisson jensen concrete}.}
        \label{fig:poisson-jensen concrete example}
    \end{figure}
    This function has a poles of multiplicities $3$ and $1$ at $0$ and $2$ respectively and a root of multiplicity $3$ at $-1$. In addition $f$ has negative jumps of height $1$ at $0$ and $-2$ and a positive jump of height $3$ at $1$. Let $x\in(-1,0)$. Denote the slopes of the function $f$ by $m_{-2} = -1,\ m_{-1} = m_1 = 2,\ m_2 = -1$ and  $m_3 = -2$, the points of discontinuity of the derivative of $f$ by $c_{-1} = -1,\ c_1 = 0$ and $c_2 = 2$ and the points of discontinuity of $f$ by $k_1 = -2,\ k_2 =0,\ k_3=1$ and $k_4 = 4$. Then we can write
    \begin{equation*}\begin{aligned}
        &f(r) - f(x) \\
        = &m_1(c_1-x) + m_2(c_2-c_1) + m_3(r-c_2) + \dcm_f(k_2) + \dcm_f(k_3) + \dcm_f(k_4)\\
        = &c_1(m_1-m_2) + c_2(m_2-m_3) + m_3r-m_1x+ \dcm_f(k_2) + \dcm_f(k_3) + \dcm_f(k_4)\\
        = &m_1(r-x) - (m_1-m_2)(r-c_1) - (m_2-m_3)(r-c_2) + \dcm_f(k_2) + \dcm_f(k_3) + \dcm_f(k_4)\\
    \end{aligned}\end{equation*}
    and similarly
    \begin{equation*}
        f(x) - f(-r) = m_{-1}(r+x)+(m_{-2}-m_{-1})(r+c_{-1}) + \dcm_f(k_1).
    \end{equation*}
    By multiplying the above equations by $(r+x)$ and $(r-x)$ respectively and subtracting we obtain
    \begin{equation*}\begin{aligned}
        &2rf(x)\\ 
        &= r(f(r) + f(-r))+ x(f(r)-f(-r)) \\\
        & + (m_{-2}-m_{-1})(r^2-r(c_{-1}-x) - c_{-1}x) +\sum_{j=1}^2 (m_j-m_{j+1})(r^2-r(x-c_j) - c_jx)\\
        & +r(\dcm_f(k_1)-\dcm_f(k_2)-\dcm_f(k_3) - \dcm_f(k_4)) - x(\dcm_f(k_1) + \dcm_f(k_2) + \dcm_f(k_3) + \dcm_f(k_4))\\
        &= r(f(r) + f(-r))+ x(f(r)-f(-r)) \\\
        & + \omega_f(c_{-1})(r^2-r(c_{-1}-x) - c_{-1}x) +\sum_{j=1}^2 \omega_f(c_j)(r^2-r(x-c_j) - c_jx)\\
        & + r(\dcm_f(k_1) - \dcm_f(k_2) - \dcm_f(k_3) - \dcm_f(k_4)) - x(\dcm_f(k_1)+\dcm_f(k_2) + \dcm_f(k_3) + \dcm_f(k_4)).
    \end{aligned}\end{equation*}
    By the definition of multiplicity of a root and a pole and the height of a jump we can see that
    \begin{equation*}
        \begin{aligned}
            &f(x) = \frac{1}{2}(f(r) + f(-r)) + \frac{x}{2r}(f(r)-f(-r))\\
            & - \frac{1}{2r}\tau_f(c_{-1})(r^2-r|c_{-1}-x| - c_{-1}x) + \frac{1}{2r}\sum_{j=1}^2 \tau_f(c_j)(r^2-r|c_j-x| - c_jx)\\
            & - \frac{1}{2}(\dcmg_f(k_3) + \dcmg_f(k_4)) + \frac{1}{2}(\dcmg_f(k_1) + \dcmg_f(k_2))  \\
            & - \frac{x}{2r}(\dcmg_f(k_1)  + \dcmg_f(k_3) + \dcmg_f(k_4))+\frac{x}{2r}\dcmg_f(k_2).
        \end{aligned}
    \end{equation*}
    Finally with concrete values we obtain
    \begin{equation*}
        \begin{aligned}
            f(x) =& \frac{1}{2}(-1 + 2) + \frac{x}{8}(-1-2)\\
            & - \frac{1}{8}3(16-4|-1-x| + x) + \frac{1}{8} 3(16-4|x|) + \frac{1}{8}(16-4|2-x| - 2x)\\
            &+ \frac{1}{2}(1 + 1 - 3 - 1) - \frac{x}{8}(1 - 1 + 3 + 1)\\
            =& 2x + 2,
        \end{aligned}
    \end{equation*}
    when $-1<x<0$.
\end{example}
\begin{proof}[Proof of Theorem \ref{theorem: poisson jensen extended}]
    Let $x\in(-r,r)$. Define an increasing sequence $(c_j)$, $j = -p,\ldots,q$ in $(-r,r)$ in the following way. Let $c_0 = x$, and let the other points in this sequence be the points in $(-r,r)$, at which the derivative of $f$ does not exist, i.e. the roots and poles of $f$. Then define 
    \begin{equation*}
        m_{j-1} = f'(c_j-) 
    \end{equation*}
    for $j = -p,\ldots,0$ and 
    \begin{equation*}
        m_{j+1} = f'(c_j+) 
    \end{equation*}
    for $j = 0,\ldots,q$. Let $(k_j), j = 1,\ldots,s$ be an increasing sequence of the discontinuities of $f$ in $[-r,r]$ including the possible discontinuity at $r$ (resp. at $-r$) only if $f$ is left-discontinuous at $r$ (resp. right-discontinuous at $-r$). Define
    \begin{equation*}
        K_j = f(k_j+)-f(k_j-)
    \end{equation*}
    for $j = 1,\ldots,s$. Also define the index sets 
    \begin{equation*}
        J^+ = \{j\in\{1,\ldots,s\}:(x\leq k_j\leq r)\land (k_j = x \implies f \text{ is left-discontinuous at }x)\}
    \end{equation*}
    and
    \begin{equation*}
        J^- = \{j\in\{1,\ldots,s\}:(-r\leq k_j\leq x) \land (k_j = x \implies f \text{ is right-discontinuous at }x)\}.
    \end{equation*}
    \begin{figure}
        \centering
        \begin{tikzpicture}
\def\xmin{-6}
\def\xmax{10}

\def\r{6}

\def\kix{-\r}
\def\kiyo{2.5}
\def\kiyv{1}

\def\cpx{-\r+0.5}
\def\cpy{2.5}

\def\kiix{-\r+1.5}
\def\kiiyo{2.5}
\def\kiiyv{0.5}

\def\negleftcutx{-\r+2}
\def\negleftcuty{1.5}

\def\negrightcutx{-\r+3}
\def\negrightcuty{2.5}

\def\xx{-2}

\def\kijx{\xx}
\def\kijyo{1.5}
\def\kijyv{3.5}

\def\kjx{-1}
\def\kjyo{3}
\def\kjyv{0.5}

\def\cix{0.5}
\def\ciy{1.5}

\def\ciix{1.5}
\def\ciiy{2}

\def\kjix{2}
\def\kjiyo{1}
\def\kjiyv{2}

\def\posleftcutx{\r-3}
\def\posleftcuty{1}

\def\posrightcutx{\r-2}
\def\posrightcuty{1.5}

\def\ksx{\r-1}
\def\ksyo{4}
\def\ksyv{1}

\def\cqx{\ksx}

\def\ksix{\r}
\def\ksiyo{2}
\def\ksiyv{0.5}

\tkzInit[xmax=\r,ymax=10,xmin=-\r,ymin=-10]
   \begin{scope}
    \end{scope}
   \tkzDrawX[noticks = true, label = ]

\tkzHTicks[mark = |]{-\r,\cpx,\kiix,\r,\xx,\kjx,\cix,\ciix,\kjix,\ksx,\ksix,0};
\node[circle, label = {-90:$-r$}] at (-\r,0) {};
\node[circle, label = {90:$k_{1}$}] at (\kix,0) {};
\node[circle, label = {-90:$c_{-p}$}] at (\cpx,0) {};
\node[circle, label = {-90:$k_{2}$}] at (\kiix,0) {};
\node[circle, label = {-90:$r$}] at (\r,0) {};
\node[circle, label = {90:$k_{j-1}$}] at (\kijx,0) {};
\node[circle, label = {-90:$x$}] at (\xx,0) {};
\node[circle, label = {-90:$k_{j}$}] at (\kjx,0) {};
\node[circle, label = {-90:$c_1$}] at (\cix,0) {};
\node[circle, label = {-90:$c_2$}] at (\ciix,0) {};
\node[circle, label = {-90:$k_{j+1}$}] at (\kjix,0) {};
\node[circle, label = {90:$k_{s}$}] at (\ksx,0) {};
\node[circle, label = {-90:$c_{q}$}] at (\cqx,0) {};
\node[circle, label = {-90:$0$}] at (0,0) {};

\draw (\kix,\kiyo) -- (\cpx,\cpy) node[midway,above]{$m_{-p-1}$};
\draw[dashed] (\kix,\kiyv) -- (\kix,\kiyo) node[midway,left]{$K_1$};
\filldraw[color = black, fill = white] (\kix,\kiyo) circle (1pt);
\filldraw (\kix,\kiyv) circle (1pt);

\draw (\cpx,\cpy) -- (\kiix,\kiiyv) node[midway,left]{$m_{-p}$};

\draw (\kiix,\kiiyo) -- (\negleftcutx,\negleftcuty) node[midway,right]{$m_{-p}$};

\draw[dashed] (\kiix,\kiiyv) -- (\kiix,\kiiyo) node[midway,right]{$K_2$};
\filldraw[color = black, fill = white] (\kiix,\kiiyo) circle (1pt);
\filldraw (\kiix,\kiiyv) circle (1pt);

\draw (\negrightcutx,\negrightcuty) -- (\kijx,\kijyv) node[midway,left]{$m_{-1}$};

\draw[dashed] (\kijx,\kijyv) -- (\kijx,\kijyo) node[midway,right]{$K_{j-1}$};
\filldraw[color = black, fill = white] (\kijx,\kijyv) circle (1pt);
\filldraw (\kijx,\kijyo) circle (1pt);

\draw (\kijx,\kijyo) -- (\kjx,\kjyv) node[midway,left]{$m_{1}$};

\draw[dashed] (\kjx,\kjyv) -- (\kjx,\kjyo) node[midway,right]{$K_j$};
\filldraw[color = black, fill = white] (\kjx,\kjyv) circle (1pt);
\filldraw (\kjx,\kjyo) circle (1pt);

\draw (\kjx,\kjyo) -- (\cix,\ciy) node[midway,right]{$m_{1}$};

\draw (\cix,\ciy) -- (\ciix,\ciiy) node[midway,above]{$m_{2}$};

\draw (\ciix,\ciiy) -- (\kjix,\kjiyv) node[midway,above]{$m_{3}$};

\draw[dashed] (\kjix,\kjiyv) -- (\kjix,\kjiyo) node[midway,left]{$K_{j+1}$};
\filldraw[color = black, fill = white] (\kjix,\kjiyv) circle (1pt);
\filldraw (\kjix,\kjiyo) circle (1pt);

\draw (\kjix,\kjiyo) -- (\posleftcutx,\posleftcuty) node[midway,above]{$m_{3}$};

\draw (\posrightcutx,\posrightcuty) -- (\ksx,\ksyv) node[midway,above]{$m_{q}$};

\draw[dashed] (\ksx,\ksyv) -- (\ksx,\ksyo) node[midway,left]{$K_{s}$};
\filldraw[color = black, fill = white] (\ksx,\ksyv) circle (1pt);
\filldraw (\ksx,\ksyo) circle (1pt);

\draw (\ksx,\ksyo) -- (\ksix,\ksiyv) node[midway,right]{$m_{q+1}$};

\filldraw (\ksix,\ksiyv) circle (1pt);

\node[circle] at (-\r+2.5,2.5) {$\cdots$};

\node[circle] at (\r-2.5,1) {$\cdots$};






\end{tikzpicture}
        \caption{Notation in the proof of the Poisson-Jensen formula.}
        \label{fig:poisson-jensen}
    \end{figure}
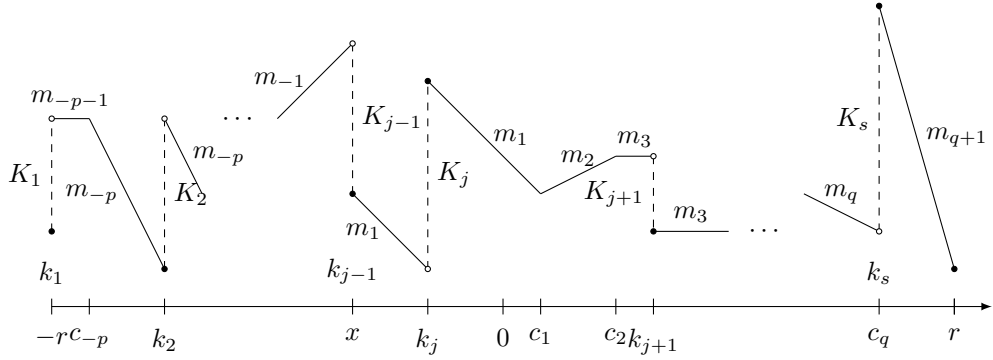
    See Figure \ref{fig:poisson-jensen} for the notation. In Figure \ref{fig:poisson-jensen} $J^- = \{1,\ldots,j-1\}$ and $J^+ = \{j,\ldots,s\}$. A geometric observation implies that
    \begin{equation*}\begin{aligned}
        &f(r) - f(x)\\
        &= m_1(c_1-x) + m_2(c_2-c_1) + \cdots + m_q(c_q-c_{q-1})+ m_{q+1}(r-c_q) + \sum_{j\in J^+}K_j\\
        &= m_1(r-x) -\sum_{j = 1}^q(m_j-m_{j+1})(r-c_j) + \sum_{j\in J^+}K_j.
    \end{aligned}\end{equation*}
    Similarly
    \begin{equation*}
        f(x) - f(-r) = m_{-1}(r+x) +\sum_{j = 1}^p(m_{-j-1}-m_{-j})(r+c_{-j}) + \sum_{j\in J^-}K_j.
    \end{equation*}
    When multiplying the equations above by $(r+x)$ and $(r-x)$, respectively, and subtracting we obtain
    \begin{equation*}\begin{aligned}
        2rf(x) 
        =& r(f(r)+f(-r)) + x(f(r)-f(-r)) + (m_{-1}-m_1)(r^2-x^2)\\
         & +\sum_{j=1}^p(m_{-j-1}-m_j)(r^2-(x-c_{-j})r-c_{-j}x)\\
         & +\sum_{j=1}^q(m_{j}-m_{j+1})(r^2-(c_j-x)r-c_{j}x)\\
         & + \sum_{j\in J^-}K_j(r-x) - \sum_{j\in J^+}K_j(r+x) \\
        =& r(f(r)+f(-r)) + x(f(r)-f(-r))\\
         & +\sum_{j=-p}^q-\omega_f(c_j)(r^2-|c_j-x|r-c_jx)\\
         & + x\sum_{j=1}^s -\dcm_f(k_j) - r\left(\sum_{j\in J^-}-\dcm_f(k_j)-\sum_{j\in J^+}-\dcm_f(k_j)\right).
    \end{aligned}\end{equation*}
    We can see that
    \begin{equation*}
        \sum_{j\in J^-}-\dcm_f(k_j)-\sum_{j\in J^+}-\dcm_f(k_j) = \sum_{-r\leq k_j\leq x}-\dcm_f(k_j)-\sum_{x\leq k_j\leq r}-\dcm_f(k_j) + A_f(x),
    \end{equation*}
    where 
    \begin{equation*}
        A_f(x) = 
        \begin{cases}
            0, &\text{if }f \text{ is continuous at }x,\\
           \dcm_f(x), &\text{if }f \text{ is left-discontinuous at }x,\\
           -\dcm_f(x), &\text{if }f \text{ is right-discontinuous at }x.
        \end{cases}
    \end{equation*}
    We also have
    \begin{equation*}\begin{aligned}
        &\sum_{-r\leq k_j\leq x}-\dcm_f(k_j)-\sum_{x\leq k_j\leq r}-\dcm_f(k_j) \\
        &= \sum_{-r\leq k_j\leq \min\{0,x\}}\!\!\!\!\!\!\!\!\!\!\!\!-\dcm_f(k_j) + \sum_{0\leq k_j\leq x}\!\!\!-\dcm_f(k_j)\\
        &-\left(\sum_{\max\{0,x\}\leq k_j\leq r}\!\!\!\!\!\!\!\!\!\!\!\!-\dcm_f(k_j) + \sum_{x\leq k_j\leq 0}\!\!\!-\dcm_f(k_j)\right) + B_f(x)\\
        &=\sum_{-r\leq \alpha_\mu\leq \min\{0,x\}}\!\!\!\!\!\!\!\!\!\!\!\!\dcmg_f(\alpha_\mu)
        -\sum_{-r\leq \beta_\nu\leq \min\{0,x\}}\!\!\!\!\!\!\!\!\!\!\!\!\dcmg_f(\beta_\nu)
        -\sum_{0\leq \alpha_\mu\leq x}\!\!\!\dcmg_f(\alpha_\mu)
        +\sum_{0\leq \beta_\nu\leq x}\!\!\!\dcmg_f(\beta_\nu)\\
        &-\left(
        -\!\!\!\!\!\sum_{\max\{0,x\}\leq \alpha_\mu\leq r}\!\!\!\!\!\!\!\!\!\!\!\!\dcmg_f(\alpha_\mu)
        +\!\!\!\!\!\sum_{\max\{0,x\}\leq \beta_\nu\leq r}\!\!\!\!\!\!\!\!\!\!\!\!\dcmg_f(\beta_\nu) 
        +\sum_{x\leq \alpha_\mu\leq 0}\!\!\!\dcmg_f(\alpha_\mu)
        -\sum_{x\leq \beta_\nu\leq 0}\!\!\!\dcmg_f(\beta_\nu)\right) + B_f(x)\\
        &=\sum_{-r\leq \alpha_\mu\leq \min\{0,x\}}\!\!\!\!\!\!\!\!\!\!\!\!\dcmg_f(\alpha_\mu) 
        + \sum_{\max\{0,x\}\leq \alpha_\mu\leq r}\!\!\!\!\!\!\!\!\!\!\!\!\dcmg_f(\alpha_\mu)
        - \sum_{x\leq \alpha_\mu\leq 0}\!\!\!\dcmg_f(\alpha_\mu)
        - \sum_{0\leq \alpha_\mu\leq x}\!\!\!\dcmg_f(\alpha_\mu)\\
        &-\sum_{-r\leq \beta_\nu\leq \min\{0,x\}}\!\!\!\!\!\!\!\!\!\!\!\!\dcmg_f(\beta_\nu)
        - \sum_{\max\{0,x\}\leq \beta_\nu\leq r}\!\!\!\!\!\!\!\!\!\!\!\!\dcmg_f(\beta_\nu) 
        + \sum_{0\leq \beta_\nu\leq x}\!\!\!\dcmg_f(\beta_\nu)
        + \sum_{x\leq \beta_\nu\leq 0}\!\!\!\dcmg_f(\beta_\nu) + B_f(x),
    \end{aligned}\end{equation*}
    where 
    \begin{equation*}
        B_f(x) = 
        \begin{cases}
            0, &\text{if }x = 0,\\
           -\dcm_f(0), &\text{if }x< 0,\\
           \dcm_f(0), &\text{if }x > 0.
        \end{cases}
    \end{equation*}
    We also have
    \begin{equation*}
        \begin{aligned}
            &\sum_{j=1}^s -\dcm_f(k_j)\\
            &= \sum_{-r\leq\alpha_\mu\leq 0}\dcmg(\alpha_\mu) - \sum_{-r\leq\beta_\nu\leq 0}\dcmg(\beta_\nu) \\
            &- \sum_{0\leq\alpha_\mu\leq r}\dcmg(\alpha_\mu) + \sum_{0\leq\beta_\nu\leq r}\dcmg(\beta_\nu) - \dcm_f(0)\\
            &= \sum_{-r\leq\alpha_\mu\leq 0}\dcmg(\alpha_\mu)+ \sum_{0\leq\beta_\nu\leq r}\dcmg(\beta_\nu)\\
            &-\left(\sum_{-r\leq\beta_\nu\leq 0}\dcmg(\beta_\nu) + \sum_{0\leq\alpha_\mu\leq r}\dcmg(\alpha_\mu)\right) - \dcm_f(0).
        \end{aligned}
    \end{equation*}

    By combining the above equalities we obtain \eqref{equation: poisson jensen statement}
    and in the case $x=0$ we have
    \begin{equation*}
        \begin{aligned}
            f(0)
            & = \frac{1}{2}(f(r)+f(-r))\\
            & - \frac{1}{2r}\sum_{|a_\mu|<r}\tau_f(a_\mu)(r^2 - |a_\mu|r)\\
            & + \frac{1}{2r}\sum_{|b_\nu|<r}\tau_f(b_\nu)(r^2 - |b_\nu|r)\\
            & - \frac{1}{2}\left(\sum_{-r\leq \alpha_\mu\leq 0}\dcmg_f(\alpha_\mu) 
              + \sum_{0\leq \alpha_\mu\leq r}\dcmg_f(\alpha_\mu)
              - 2\sum_{\alpha_\mu = 0}\dcmg_f(\alpha_\mu)\right)\\
            & +\frac{1}{2}\left(\sum_{-r\leq \beta_\nu\leq 0}\dcmg_f(\beta_\nu)
              + \sum_{0\leq \beta_\nu\leq r}\dcmg_f(\beta_\nu) 
              - 2\sum_{\beta_\nu = 0}\dcmg_f(\beta_\nu)\right)\\
            & -\frac{1}{2}A_f(0)\\
            & = \frac{1}{2}(f(r)+f(-r))
             - \frac{1}{2}\sum_{|a_\mu|<r}\tau_f(a_\mu)(r - |a_\mu|)
             + \frac{1}{2}\sum_{|b_\nu|<r}\tau_f(b_\nu)(r - |b_\nu|)\\
            & - \frac{1}{2}\sum_{|\alpha_\mu|\leq r}\dcmg_f(\alpha_\mu)
             +\frac{1}{2}\sum_{|\beta_\nu|\leq r}\dcmg_f(\beta_\nu)\\
            & -\frac{1}{2}\left(A_f(0)-\sum_{\alpha_\mu = 0}\dcmg_f(\alpha_\mu)+\sum_{\beta_\nu = 0}\dcmg_f(\beta_\nu)\right).
        \end{aligned}
    \end{equation*}
    By the definition of the positive and negative jump at $0$ we can see that
    \begin{equation*}
        A_f(0)-\sum_{\alpha_\mu = 0}\dcmg_f(\alpha_\mu)+\sum_{\beta_\nu = 0}\dcmg_f(\beta_\nu) = 0.
    \end{equation*}
\end{proof}

For piecewise linear functions with discontinuities we will define the proximity function and counting function in the same way as for tropical meromorphic functions. The proximity function and the counting function are well defined and the properties
\begin{align*}
    V(r,f^{\otimes\alpha}) &= \alpha V(r,f)\\
    V(r,f\otimes g)&\leq V(r,f) + V(r,g)\\
    V(r,f\oplus g)&\leq V(r,f) + V(r,g)
\end{align*}
remain true for both $V(r,f)=m(r,f)$ and $V(r,f) = N(r,f)$. However, the proximity function is not continuous for noncontinuous functions. 

We will define the jump counting function as 
\begin{equation*}
    \dcn(r,f) = \frac{1}{2}\sum_{j = 1}^n\dcmg_f(\beta_j),
\end{equation*}
where $\beta_1,\ldots,\beta_n$ are the negative jumps of $f$ on the interval $[-r,r]$ counting the possible negative jump at $r$ (resp. at $-r$) if $f$ is left-discontinuous at $r$ (resp. right-discontinuous at $-r$). 

Then we can define the tropical Nevanlinna characteristic function for a piecewise linear function $f$ as 
\begin{equation*}
    T(r,f) = m(r,f) + N(r,f) + \dcn(r,f).
\end{equation*}

The jump counting function is a non-negative non-decreasing piecewise constant function. This means that the characteristic function will also remain a non-negative non-decreasing piecewise linear function. The jump counting function also has the following properties
\begin{align*}
    \dcn(r,f^{\otimes\alpha}) &= \alpha\dcn(r,f)\\
    \dcn(r,f\otimes g)&\leq\dcn(r,f) + \dcn(r,g)\\
    \dcn(r,f\oplus g)&\leq \dcn(r,f) + \dcn(r,g).
\end{align*}

This means that all the above properties remain true for the characteristic function. The only basic properties of the characteristic function that do not always hold with this definition are continuity and convexity. That is because as a piecewise constant function, the jump counting function is not always continuous or convex.  

The property
\begin{equation*}
    \dcn(r,f\oplus g)\leq \dcn(r,f) \oplus \dcn(r,g)
\end{equation*}
also does not hold. This can be seen by considering the functions
\begin{equation*}
    f(x) = \begin{cases}
        2, &x\leq 1\\
        3-2x, &x > 1
    \end{cases}
    \quad\text{ and }\quad g(x) = \begin{cases}
        \frac{1}{2}, &x \leq 2\\
        0, &x > 2.
    \end{cases}
\end{equation*}
In this case $\dcn(r,f) = 1$, $\dcn(r,g) = \frac{1}{2}$ and $\dcn(r,f\oplus g) = \frac{3}{2}$ for all $r>2$. 

Next we will move on to the second main theorem. The tropical second main theorem was first proven by Laine and Tohge. 
\begin{aptheorem}[\cite{lainet:11}]\label{theorem: second main theorem old}
    Suppose that is $f$ a non-constant tropical meromorphic function of hyper-order $\rho_2<1$, and take $0<\delta<1-\rho_2$. If $q\geq 1$ distinct values $a_1,\ldots,a_q\in\R$ satisfy 
    \begin{equation}\label{equation: laine tohge second main theorem condition about poles and target values}
        \max\{a_1,\ldots,a_q\}<\inf\{f(b):\omega_f(b)<0\}
    \end{equation}
    and 
    \begin{equation}\label{equation: old second main thm droppable assumption}
        \inf\{f(b):\omega_f(b)>0\} > -\infty.
    \end{equation}
    Then
    \begin{equation} \label{equation: old second main theorem result}
        \begin{aligned}
            qT(r,f)\leq \sum_{j=1}^qN\left(r,\ofrac{1_\circ}{f\oplus a_j}\right)+o\left(\frac{T(r,f)}{r^\delta}\right)
        \end{aligned}
    \end{equation}
    outside an exceptional set of finite logarithmic measure. 
\end{aptheorem}
Korhonen and Tohge \cite{korhonent:16AM} improved this result  by dropping the assumption \eqref{equation: old second main thm droppable assumption} and by noting that \eqref{equation: old second main theorem result} can be turned into an equality. In this paper we will show that the growth condition can be dropped entirely and in addition, we shall generalize the second main theorem for piecewise linear functions with discontinuities and for a class of piecewise linear target values instead of constant targets. 

To prove the second main theorem, we will need to prove a version of the first main theorem. Laine and Tohge \cite{lainetohge:19} proved the first main theorem for tropical meromorphic functions on a finite interval. With a similar proof we can prove the same result for piecewise linear functions on the whole real line with piecewise linear targets.
\begin{aptheorem}\label{theorem: first main theorem new}
    Let $f$ and $a$ be piecewise linear functions. Then
\begin{equation}\begin{aligned}\label{equation: first main thm statement}
    T\left(r,\ofrac{1_\circ}{f\oplus a}\right) &= T(r, f) - N(r,f) + N(r,f\oplus a) \\
    &- \dcn(r,f) + \dcn(r,f\oplus a) - f(0)\oplus a(0) + \varepsilon_f(r,a),
\end{aligned}\end{equation}
for some quantity $\varepsilon_f(r,a)$ satisfying $0\leq \varepsilon_f(r,a)\leq m(r,a)$ for all $r\geq 0$.
\end{aptheorem}
\begin{proof}
    By the Jensen formula we have
    \begin{equation*}
    \begin{aligned}
        &T\left(r,\ofrac{1_\circ}{f\oplus a}\right) -\left( T(r, f) - N(r,f) + N(r,f\oplus a) - \dcn(r,f) + \dcn(r,f\oplus a)\right)\\
        &=T\left(r,f\oplus a\right) -\left( m(r, f)  + N(r,f\oplus a) + \dcn(r,f\oplus a)\right) + f(0)\oplus a(0)\\
        &=m\left(r,f\oplus a\right) - m(r,f) + f(0)\oplus a(0).
    \end{aligned}
    \end{equation*}
    By the properties of the proximity function we obtain
    \begin{equation*}
        m\left(r,f\oplus a\right) - m(r,f)\leq m(r,f) + m(r,a)-m(r,f) = m(r,a)
    \end{equation*}
    and
    \begin{equation*}
        m\left(r,f\oplus a\right) - m(r,f)\geq m(r,f) - m(r,f) = 0.
    \end{equation*}
\end{proof}
If $m(r,a) = o(T(r,f))$, where $r$ approaches infinity, then \eqref{equation: first main thm statement} can be written in the form
\begin{equation*}
    T\left(r,\ofrac{1_\circ}{f\oplus a}\right) = T(r, f) - N(r,f) + N(r,f\oplus a)- \dcn(r,f) + \dcn(r,f\oplus a) + o(T(r,f)),
\end{equation*}
where $r$ approaches infinity. From now on in this paper, unless specified otherwise, whenever we use asymptotic notation such as $g(r) = o(T(r,f))$ or $g(r) = O(1)$ it is implied that $r$ approaches infinity without any exceptional set.

On the other hand, if $a\in\R$ is a constant we can see that
\begin{equation*}
    \ofrac{1_\circ}{f(x)\oplus a} = -\max\{f(x),a\} = \min\{-f(x),-a\} \leq -a,
\end{equation*}
which means that
\begin{equation}\label{equation: m with constant is small}
    0\leq m\left(r, \ofrac{1_\circ}{f\oplus a}\right)\leq (-a)^+.
\end{equation}
Now \eqref{equation: first main thm statement} can be written in the form
\begin{equation*}\begin{aligned}
    T(r, f) =& N\left(r,\ofrac{1_\circ}{f\oplus a}\right) + \dcn\left(r,\ofrac{1_\circ}{f\oplus a}\right) \\
    &+ N(r,f) - N(r,f\oplus a) \\
    &+ \dcn(r,f) - \dcn(r,f\oplus a) + f(0)\oplus a(0) - \varepsilon_f(r,a),
\end{aligned}\end{equation*}
where $0\leq \varepsilon_f(a,r)\leq a^+ + (-a)^+ = |a|$ for all $a\in\R$ and $r\geq 0$. This can be seen as the second main theorem for piecewise linear functions with constant target values. Motivated by this we introduce the following version of the second main theorem for a class of piecewise linear target values. 

\begin{aptheorem}\label{theorem: second main theorem}
    Let $f$ and $a$ be tropical meromorphic functions such that $m(r,a) = o(T(r,f))$. If
    \begin{equation}\label{equation: second main theorem 1/f+a}
        m\left(r,\ofrac{1_\circ}{f\oplus a}\right) = o(T(r,f))
    \end{equation}
    then
\begin{equation}\begin{aligned}\label{equation: regular second main theorem statement}
    T(r,f)&= N\left(r,\ofrac{1_\circ}{f\oplus a}\right) + \dcn\left(r,\ofrac{1_\circ}{f\oplus a}\right) \\
    &+ N(r,f) - N(r,f\oplus a)\\
    &+ \dcn(r,f) - \dcn(r,f\oplus a)  + o(T(r,f)).
\end{aligned}\end{equation}
\end{aptheorem}
    Note that by \eqref{equation: m with constant is small} the condition \eqref{equation: second main theorem 1/f+a} is always satisfied for constant values $a\in\R$. In fact, the error term $o(T(r,f))$ becomes $O(1)$ for constant values $a\in\R$.
    
    The condition \eqref{equation: second main theorem 1/f+a} is implied by $m\left(r,\ofrac{1_\circ}{f}\right)=o(T(r,f))$ or $m\left(r,\ofrac{1_\circ}{a}\right)=o(T(r,f))$. However, the converse is not true. For example if $f(x) = -x$ and $a(x) = -\max\{-x,0\}$, then 
    $$
    T(r,f) = m\left(r,\ofrac{1_\circ}{f}\right) = \frac{1}{2}r = m\left(r,\ofrac{1_\circ}{a}\right),
    $$
    but
    \begin{equation*}
        m\left(r,\ofrac{1_\circ}{f\oplus a}\right) = 0.
    \end{equation*}

Theorem \ref{theorem: second main theorem} works also if $a\equiv 0_\circ = -\infty$ and $m\left(r,\ofrac{1_\circ}{f\oplus a}\right) = m\left(r,\ofrac{1_\circ}{f}\right) = o(T(r,f))$. This can be seen by considering the Jensen formula
\begin{equation*}
    N\left(r,\ofrac{1_\circ}{f}\right) + m\left(r,\ofrac{1_\circ}{f}\right) + \dcn\left(r,\ofrac{1_\circ}{f}\right) =  N(r,f) + m(r,f) + \dcn(r,f) + O(1)
\end{equation*}
and the fact that $N(r,f\oplus 0_\circ) \equiv N(r,f)$ and $J(r,f\oplus 0_\circ) \equiv J(r,f)$.

If $f$ is a tropical meromorphic function and $a\in\R$ is a constant such that \eqref{equation: laine tohge second main theorem condition about poles and target values} holds, then we can see that $N(r,f\oplus a) \equiv N(r,f)$ and therefore
\begin{equation*}
    T(r,f) = N\left(r,\ofrac{1_\circ}{f\oplus a}\right) + O(1).
\end{equation*}
This means that Theorem \ref{theorem: second main theorem} implies Theorem \ref{theorem: second main theorem old}. In the following result we will introduce conditions similar to \eqref{equation: laine tohge second main theorem condition about poles and target values} with piecewise linear targets.
\begin{apcorollary}\label{corollary: second main thm regular corollary}
    Let $f$ and $a$ be piecewise linear functions such that $T(r,a) = o(T(r,f))$. 
    Then
    \begin{enumerate}[(i)]
        \item if $\max\{a(x+),a(x-)\}<\min\{f(x+),f(x-)\}$ for all $x\in \{x\in\R:\omega_f(x)<0\}$ we have 
        \begin{equation}\label{equation: regular second main theorem corollary eq 1}
            T(r,f) = N\left(r,\ofrac{1_\circ}{f\oplus a}\right)+\dcn\left(r,\ofrac{1_\circ}{f\oplus a}\right) + \dcn(r,f) - \dcn(r,f\oplus a) + o(T(r,f));
        \end{equation}
        \item if $\max\{a(x+),a(x-)\}<\min\{f(x+),f(x-)\}$ for all $x\in \{x\in\R:x\dcm_f(x)<0\}$ we have 
        \begin{equation*}
            T(r,f) = N\left(r,\ofrac{1_\circ}{f\oplus a}\right)+\dcn\left(r,\ofrac{1_\circ}{f\oplus a}\right) + N(r,f) - N(r,f\oplus a) + o(T(r,f));
        \end{equation*}
        \item if $\min\{a(x+),a(x-)\}>\max\{f(x+),f(x-)\}$ for all $x\in \{x\in\R:\omega_f(x)<0\}$ we have 
        \begin{equation}\label{equation: regular second main theorem corollary eq 2}
            m(r,f) = N\left(r,\ofrac{1_\circ}{f\oplus a}\right) +\dcn\left(r,\ofrac{1_\circ}{f\oplus a}\right) + \dcn(r,f) - \dcn(r,f\oplus a) + o(T(r,f));
        \end{equation}
        \item if $\min\{a(x+),a(x-)\}>\max\{f(x+),f(x-)\}$ for all $x\in \{x\in\R:x\dcm_f(x)<0\}$ we have 
        \begin{equation*}
            m(r,f) = N\left(r,\ofrac{1_\circ}{f\oplus a}\right) +\dcn\left(r,\ofrac{1_\circ}{f\oplus a}\right) + N(r,f) - N(r,f\oplus a) + o(T(r,f)).
        \end{equation*}
\end{enumerate}
\end{apcorollary}
\begin{proof}
    First we attain
    \begin{equation*}
        \ofrac{1_\circ}{f(x)\oplus a(x)} = -\max\{f(x),a(x)\} = \min\{-f(x),-a(x)\} \leq -a(x).  
    \end{equation*}
    Now because $T(r,a) = o(T(r,f))$ we can see that
    \begin{equation*}
        m\left(r,\ofrac{1_\circ}{f\oplus a}\right)\leq m\left(r,\ofrac{1_\circ}{a}\right) \leq T\left(r,\ofrac{1_\circ}{a}\right) = T(r,a) + O(1) = o(T(r,f)).
    \end{equation*}
    We can now utilize Theorem \ref{theorem: second main theorem} to obtain
    \begin{equation}\label{equation: second main theorem statement in corollary}
        T(r,f) = N\left(r,\ofrac{1_\circ}{f\oplus a}\right)+N(r,f) - N(r,f\oplus a)+\dcn(r,f)-\dcn(r,f\oplus a)+o(T(r,f)).
    \end{equation}
    Assume that $\max\{a(x+),a(x-)\}<\min\{f(x+),f(x-)\}$ for all $x\in \{x:\omega_f(x)<0\}$. Then we know that $N(r,f\oplus a)\geq N(r,f)$. In general it is true that $N(r,f\oplus a)\leq N(r,f) + N(r,a)$. So overall we have 
    \begin{equation}\label{equation: regular second main theorem corollary proof eq 1}
        0\geq N(r,f) - N(r,f\oplus a)\geq -N(r,a) = o(T(r,f)).
    \end{equation} 
    Now \eqref{equation: regular second main theorem corollary eq 1} follows by combining \eqref{equation: second main theorem statement in corollary} and \eqref{equation: regular second main theorem corollary proof eq 1}. Part (ii) follows from similar reasoning.

    Assume next that $\min\{a(x+),a(x-)\}>\max\{f(x+),f(x-)\}$ for all $x\in \{x:\omega_f(x)<0\}$. Then we can see that $0\leq N(r,f\oplus a) \leq N(r,a) = o(T(r,f))$. Now by subtracting $N(r,f)$ from both sides of \eqref{equation: second main theorem statement in corollary} we obtain \eqref{equation: regular second main theorem corollary eq 2}. Part (iv) follows from similar reasoning.

\end{proof}

    By the proof of Corollary \ref{corollary: second main thm regular corollary} we can see that the assumptions $m(r,a) = o(T(r,f))$ and $m\left(r,\ofrac{1_\circ}{f\oplus a}\right) = o(T(r,f))$ of Theorem \ref{theorem: second main theorem} follow from the assumption $T(r,a) = o(T(r,f))$ of Corollary \ref{corollary: second main thm regular corollary}. However the converse is not true. For example, if we choose $f(x) = -|x|+1$ and $a(x) = \max\{-|x+2|+1,-|x-2|+1,0\}$ we can see that
    \begin{equation*}
        m(r,a) = O(1)\quad \text{ and }\quad m\left(r,\ofrac{1_\circ}{f\oplus a}\right) = O(1).
    \end{equation*}
    On the other hand, $N(r,a) = 2(r-2)^+$ and $T(r,f) = r$, which means $T(r,a)\neq o(T(r,f))$. In addition, even though the assumption $a(x)<f(x)$ for all $x\in\{b:\omega_f(b)<0\}$ holds, we have
    \begin{equation*}
        N\left(r,\ofrac{1_\circ}{f\oplus a}\right) = \begin{cases}0, &0\leq r<1,\\2r-2, &1\leq r\leq 3,\\3r-5, &r>3,\end{cases}
    \end{equation*}
    which means that 
    \begin{equation*}
        N\left(r,\ofrac{1_\circ}{f\oplus a}\right) = 3T(r,f) + O(1).
    \end{equation*}
    This means that we cannot weaken the assumptions of Corollary \ref{corollary: second main thm regular corollary} (i) to the same assumptions as in Theorem \ref{theorem: second main theorem}. The same can be said about Corollary \ref{corollary: second main thm regular corollary} (iii) by considering the  function $f(x) = -|x|-2$ instead, and the same function for $a(x)$ as previously.
    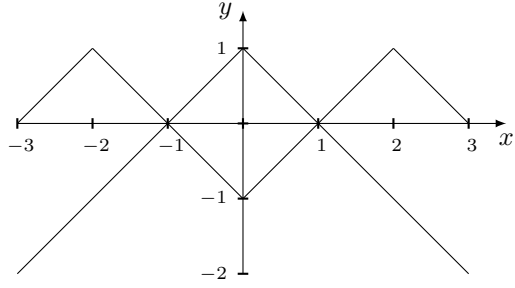
\begin{figure}
        \centering
        \begin{tikzpicture}
\def\xx{3}
\def\yy{1}
\tkzInit[xmax=\xx,ymax=1,xmin=-\xx,ymin=-2]
   \begin{scope}
    \end{scope}
   \tkzLabelX[orig=false,step=1,node font=\scriptsize]
   \tkzLabelY[orig=false,step=1,node font=\scriptsize]
   \tkzDrawX
   \tkzDrawY

\draw [domain=0:\xx] plot(\x,-\x + 1) node[right]{};
\draw [domain=-\xx:0] plot(\x,\x+1) node[right]{};

\draw [domain=0:2] plot(\x,\x-1) node[right]{};
\draw [domain=2:\xx] plot(\x,-\x+3) node[right]{};

\draw [domain=-2:0] plot(\x,-\x-1) node[right]{};
\draw [domain=-\xx:-2] plot(\x,\x+3) node[right]{};








\end{tikzpicture}    
        \caption{Functions $f(x) = -|x|+1$ and $a(x) = \max\{-|x+2|+1,-|x-2|+1,0\}$.}
        \label{fig:my_label}
    \end{figure}

\begin{apremark}
Theorem \ref{theorem: second main theorem} and Corollary \ref{corollary: second main thm regular corollary} also work on finite intervals. The proofs are identical, so they are not presented here. For more information on tropical Nevanlinna theory on finite intervals see \cite{lainetohge:19}.
\end{apremark}

\section{The second main theorem with tropical hypersurfaces}
Korhonen and Tohge  proved the tropical version of the Cartan's second main theorem. 
\begin{aptheorem}[\cite{korhonent:16AM}]
    Let $q$ and $n$ be positive integers with $q>n$, and let $\varepsilon>0$. Given $n+1$ tropical entire functions $g_0,\ldots,g_n$ without common roots, and linearly independent in Gondran-Minoux sense, let the $q+1$ tropical linear combinations $f_0\ldots,f_q$ of the $g_j$ over the semi-ring $\T$ be defined by
    \begin{equation*}
        f_k(x) = a_{0k}\otimes g_0\oplus a_{1k}\otimes g_1(x)\oplus\cdots\oplus a_{nk}\otimes g_n(x),\quad 0\leq k\leq q.
    \end{equation*}
    Let $\lambda = \operatorname{ddg}(\{f_{n+1},\ldots,f_q\})$ and
    \begin{equation*}
        L = \ofrac{f_0\otimes f_1\otimes\cdots\otimes f_n\otimes f_{n+1}\otimes\cdots\otimes f_q}{C_\circ(f_0,f_1,\ldots,f_n)}.
    \end{equation*}
    If the tropical holomorphic curve $g$ of $\R$ into $\mathbb{TP}^n$ with reduced representation $\mathbf{g}= (g_0,\ldots,g_n)$ is of hyper-order
    \begin{equation*}
        \rho_2 := \rho_2(\mathbf{g}) < 1,
    \end{equation*}
    then
    \begin{equation*}
        (q-n-\lambda)T_g(r)\leq N\left(r,\ofrac{1_\circ}{L}\right)-N(r,L) + o\left(\frac{T_g(r)}{r^{1-\rho_2-\varepsilon}}\right),
    \end{equation*}
    where $r$ approaches infinity outside an exceptional set of finite logarithmic measure. 
\end{aptheorem}
Cao and Zheng improved this result by extending it for tropical hypersurfaces and by improving the growth condition. 

\begin{aptheorem}[\cite{caozheng:18}]\label{theorem: hypersurface second main theorem old}
    Let $q$, $n$ and $d$ be positive integers such that $q > M$, where $M=\binom{n+d}{d}-1$. Let the tropical holomorphic curve $f:\R\to\mathbb{TP}^n$ be tropical algebraically nondegenerate. Assume that tropical hypersurfaces $V_{P_0},\ldots,V_{P_q}$ are defined by homogeneous tropical polynomials $P_0,\ldots,P_q$ with degrees $d_0,\ldots,d_q$, respectively, such that the least common multiple of $d_0,\ldots,d_q$ is $d$. If $\lambda = \operatorname{ddg}(\{P_{M+1}\circ f,\ldots,P_{q}\circ f\})$ and
    \begin{equation*}
        \limsup_{r\to\infty}\frac{\log T_f(r)}{r} = 0,
    \end{equation*}
    then
    \begin{equation*}
        \begin{aligned}
            &(q-M-\lambda)T_f(r)\\
            \leq &\sum_{j=0}^q\frac{1}{d_j}N\left(r,\ofrac{1_\circ}{P_j\circ f}\right) \\
            &- \frac{1}{d}N\left(r,\ofrac{1_\circ}{C_\circ(P_0^{\otimes\frac{d}{d_0}}\circ f,\ldots,P_{M}^{\otimes\frac{d}{d_{M}}}\circ f)}\right) + o(T_f(r))\\
            = &\sum_{j=M+1}^q \frac{1}{d_j}N\left(r,\ofrac{1_\circ}{P_j\circ f}\right) + o(T_f(r))\\
            \leq &(q-M)T_f(r),
        \end{aligned}
    \end{equation*}
    where $r$ approaches infinity outside an exceptional set of zero upper density measure. In the special case whenever $\lambda = 0$, 
    \begin{equation*}
        \begin{aligned}
            &(q-M)T_f(r)\\
            = &\sum_{j=0}^q \frac{1}{d_j} N\left(r,\ofrac{1_\circ}{P_j\circ f}\right) \\
            &- \frac{1}{d}N\left(r,\ofrac{1_\circ}{C_\circ(P_0^{\otimes\frac{d}{d_0}}\circ f,\ldots,P_{M}^{\otimes\frac{d}{d_{M}}}\circ f)}\right) + o(T_f(r))\\
            = &\sum_{j=M+1}^q \frac{1}{d_j}N\left(r,\ofrac{1_\circ}{P_j\circ f}\right) + o(T_f(r)),\\
        \end{aligned}
    \end{equation*}
    where $r$ approaches infinity outside an exceptional set of zero upper density measure.
\end{aptheorem}

In Theorem \ref{theorem: hypersurface second main theorem old} the homogeneous tropical polynomials have integer exponents and constant coefficients. In this paper we will consider tropical homogeneous polynomials with non-negative real number exponents and tropical meromorphic coefficients. Such a polynomial can be written in the form
\begin{equation*}
    P(x,x_0,\ldots,x_n) = \bigoplus_{\alpha_0+\alpha_1+\cdots+\alpha_n = d}a_{\alpha_0,\ldots,\alpha_n}(x)\otimes x_0^{\otimes\alpha_0}\otimes  x_1^{\otimes\alpha_1}\otimes\cdots\otimes x_n^{\otimes\alpha_n},
\end{equation*}
where $a_{\alpha_0,\ldots,\alpha_n}(x)$ are tropical meromorphic functions such that $a_{\alpha_0,\ldots,\alpha_n}(x) \not\equiv 0_\circ$ for only finitely many $(\alpha_0,\ldots,\alpha_n)\in\R_{+}^{n+1}$ such that $\alpha_0+\alpha_1+\cdots+\alpha_n = d$. If $f = [f_0:\cdots:f_n]:\R\to\mathbb{TP}^n$ is a tropical holomorphic curve then we often consider the composition of a tropical homogeneous polynomial $P$ and $f$
\begin{equation*}
    P(f)(x) = \bigoplus_{\alpha_0+\alpha_1+\cdots+\alpha_n = d}a_{\alpha_0,\ldots,\alpha_n}(x)\otimes f_0(x)^{\otimes\alpha_0}\otimes  f_1(x)^{\otimes\alpha_1}\otimes\cdots\otimes f_n(x)^{\otimes\alpha_n}.
\end{equation*}

Next we will define
\begin{equation}
    \operatorname{\psi}(P,f) := \liminf_{r\to\infty}\frac{\frac{1}{2}\left(P(f)(r) + P(f)(-r)\right)}{dT_f(r)},
\end{equation}
and
\begin{equation}
    \operatorname{\Psi}(P,f) := \limsup_{r\to\infty}\frac{\frac{1}{2}\left(P(f)(r) + P(f)(-r)\right)}{dT_f(r)},
\end{equation}
for a homogeneous tropical polynomial $P$ and a tropical holomorphic curve $f$. These values will appear in many results in this paper.

The next theorem improves the results above by dropping the growth condition and by making the inequalities tighter.
\begin{aptheorem}\label{theorem: second main theorem hypersurface new}
Let $f : \R \to \mathbb{TP}^n$ be a non-constant tropical holomorphic curve and let $P:\mathbb{TP}^n\to\R(\not\equiv 0_\circ)$ be a tropical homogeneous polynomial of degree $d$ with tropical meromorphic coefficients $a_0(x),\ldots,a_k(x)$. If for all coefficients $a_j$ of $P\circ f$ we have $N(r,a_j) = o(T_f(r))$, then 
\begin{equation*}
    \operatorname{\psi}(P,f)T_f(r) \leq \frac{1}{d}N\left(r,\frac{1_\circ}{P\circ f}\oslash\right) + o(T_f(r)) \leq \operatorname{\Psi}(P,f)T_f(r).
\end{equation*}
\end{aptheorem}
\begin{proof}
Since $f$ is a tropical holomorphic curve $P\circ f$ can have poles only at the poles of the meromorphic coefficients of $P$ and therefore $N(r,P\circ f)=o(T_f(r))$. By the tropical Jensen formula we have
    \begin{equation*}
        \begin{aligned}
    &N\left(r,\frac{1_\circ}{P\circ f}\oslash\right) \\
    &= \frac{1}{2}\left( P(f)(r)+P(f)(-r)\right)-P(f)(0)  + N(r,P\circ f)\\
    &= \frac{1}{2}\left( P(f)(r)+P(f)(-r)\right) + o(T_f(r)).
    \end{aligned}
    \end{equation*}
    We can see that
    \begin{equation*}
        \frac{1}{2}\left( P(f)(r)+P(f)(-r)\right)\leq \Psi(P,f) dT_f(r) + o(T_f(r))
    \end{equation*}
    and
    \begin{equation*}
        \frac{1}{2}\left( P(f)(r)+P(f)(-r)\right)\geq \psi(P,f) dT_f(r) + o(T_f(r)).
    \end{equation*}
    Therefore
    \begin{equation*}
        \psi(P,f) T_f(r) \leq \frac{1}{d}N\left(r,\frac{1_\circ}{P\circ f}\oslash\right)+ o(T_f(r))\leq \Psi(P,f)T_f(r).
    \end{equation*}
\end{proof}
The next lemma gives bounds for the values $\psi(P,f)$ and $\Psi(P,f)$.
\begin{aplemma}\label{lemma: liminf limsup bounds}
    Let $f : \R \to \mathbb{TP}^n$ be a non-constant tropical holomorphic curve and let $P:\mathbb{TP}^n\to\R(\not\equiv 0_\circ)$ be a homogeneous tropical polynomial with tropical meromorphic coefficients $a_1,\ldots,a_k$. If for all coefficients $a_j$ we have $m(r,a_j) = o(T_f(r))$, then
\begin{equation*}
    0\leq \psi(P,f)\leq \Psi(P,f)\leq 1.
\end{equation*}
\end{aplemma}

\begin{proof}
    Since $P$ is homogeneous, we can write $P\circ f$ in the form
    \begin{equation*}
        P(f)(x) = \max_{\alpha_0+\cdots+\alpha_n=d}\{a_{\alpha_0,\ldots,\alpha_n}(x) + \alpha_0f_0(x)+\cdots+\alpha_nf_n(x)\},
    \end{equation*}
    where $a_{\alpha_0,\ldots,\alpha_n}(x)\not\equiv 0_\circ$ for finitely many $(\alpha_0,\ldots,\alpha_n)\in\R_+^{n+1}$ such that $\alpha_0+\cdots+\alpha_n = d$. Then since at least one coefficient is non-zero we have
    \begin{equation*}
    \begin{aligned}
        P(f)(r) &\leq  \max_{\alpha_0+\cdots+\alpha_n=d}\{\alpha_0f_0(r)+\cdots+\alpha_nf_n(r)\}+\max_{\alpha_0+\cdots+\alpha_n=d}\{a_{\alpha_0,\ldots,\alpha_n}(r)\}\\
        &\leq d\|f(r)\| + o(T_f(r)).
    \end{aligned}
    \end{equation*}
    This implies that
    \begin{equation*}
    \begin{aligned}
        \Psi(P,f)
        &=\limsup_{r\to\infty}\frac{\frac{1}{2}\left(P(f)(r) + P(f)(-r)\right)}{dT_f(r)}\\
        &\leq \limsup_{r\to\infty}\frac{\frac{1}{2}d\left(\|f(r)\|+\|f(-r)\| \right) + o(T_f(r))}{dT_f(r)} = 1.
    \end{aligned}
    \end{equation*}
    Next by convexity we obtain
    \begin{equation*}
    \begin{aligned}
        &\frac{1}{2}\left(P(f)(r) + P(f)(-r)\right)\\
        &\geq \frac{1}{2}\left(\sum_{\sigma = \pm 1}\max_{\substack{\alpha_0+\cdots+\alpha_n=d\\{a_{\alpha_0,\ldots,\alpha_n}}\not\equiv0_\circ}}\{\alpha_0f_0(\sigma r)+\cdots+\alpha_nf_n(\sigma r)\} \right)\\
        &+ \frac{1}{2}\left(\sum_{\sigma = \pm 1}\min_{\substack{\alpha_0+\cdots+\alpha_n=d\\{a_{\alpha_0,\ldots,\alpha_n}}\not\equiv0_\circ}}\{a_{\alpha_0,\ldots,\alpha_n}(\sigma r)\}\right)\\
        &\geq \max_{\substack{\alpha_0+\cdots+\alpha_n=d\\{a_{\alpha_0,\ldots,\alpha_n}}\not\equiv0_\circ}}\{\alpha_0f_0(0)+\cdots+\alpha_nf_n(0)\}\\
        &+ \frac{1}{2}\left(\sum_{\sigma = \pm 1}\min_{\substack{\alpha_0+\cdots+\alpha_n=d\\{a_{\alpha_0,\ldots,\alpha_n}}\not\equiv0_\circ}}\{a_{\alpha_0,\ldots,\alpha_n}(\sigma r)\}\right)\\
        &= \max_{\substack{\alpha_0+\cdots+\alpha_n=d\\{a_{\alpha_0,\ldots,\alpha_n}}\not\equiv0_\circ}}\{\alpha_0f_0(0)+\cdots+\alpha_nf_n(0)\} + o(T_f(r)).
    \end{aligned}
    \end{equation*}
    Then 
    \begin{equation*}
    \begin{aligned}
        &\psi(P,f) = \liminf_{r\to\infty}\frac{\frac{1}{2}\left(P(f)(r) + P(f)(-r)\right)}{dT_f(r)}\\
        &\geq \liminf_{r\to\infty}\frac{\max\limits_{\alpha_0+\cdots+\alpha_n=d}^{a_{\alpha_0,\ldots,\alpha_n}\not\equiv0_\circ}\{\alpha_0f_0(0)+\cdots+\alpha_nf_n(0)\} + o(T_f(r))}{dT_f(r)} = 0.
    \end{aligned}
    \end{equation*}
\end{proof}
If we do not assume $m(r,a_j) = o(T_f(r))$ in Lemma \ref{lemma: liminf limsup bounds} the upper bound can be arbitrarily large, but the lower bound stays at $0$. 

Whenever $\psi(P,f)=\Psi(P,f)$, we obtain an equality in Theorem \ref{theorem: second main theorem hypersurface new}. In fact, it is possible that $\psi(P,f)=\Psi(P,f)=t$ for any $t\in[0,1]$. That can be seen by Lemma~\ref{lemma: for any t there exists curve and polynomial}.

    It is also possible that $\psi(P,f)<\Psi(P,f)$. To see that let us first define two tropical entire functions $f$ and $g$. Set $f(x) = g(x) = 1$ for $x\in(-\infty,0]$. Next set $g(x) = x+1$ on the interval $(0,2]$, $f(x) = 1$ on the interval $(0,1]$ and $f(x) = 2x-1$ on the interval $(1,3]$. Now we can see that $\frac{f(0)}{g(0)}=1$,$\frac{f(1)}{g(1)}=\frac{1}{2}$ and $\frac{f(2)}{g(2)}=1$. By defining $g(x) = 7x-11$ on the interval $(2,4]$, we again see that $\frac{f(3)}{g(3)}=\frac{1}{2}$. We can continue this pattern on the whole real line to obtain functions $f$ and $g$ such that $\frac{f(n)}{g(n)} = 1$ for all even $n$ and $\frac{f(n)}{g(n)} =\frac{1}{2}$ for all odd $n$ and also $g(x)\geq f(x)$ for all $x\in\R$. Since $f$ and $g$ are tropical entire functions we can define a tropical holomorphic curve $h(x) = (f(x),g(x))$ and a homogeneous tropical polynomial $P(x,y) = x\oplus (0_\circ\otimes y) = x$. Since $g(x)\geq f(x)$ for all $x\in\R$ we have $T_h(r) = \frac{1}{2}g(r) + 1$. Now we can see that
    \begin{equation*}
        \psi(P,h) = \liminf_{r\to\infty}\frac{\frac{1}{2}\left(P(h)(r) + P(h)(-r)\right)}{T_h(r)} = \liminf_{r\to\infty}\frac{f(r) + 1}{g(r) + 1} = \frac{1}{2}
    \end{equation*}
    and similarly
    \begin{equation*}
        \Psi(P,h) = \limsup_{r\to\infty}\frac{\frac{1}{2}\left(P(h)(r) + P(h)(-r)\right)}{T_h(r)} = \limsup_{r\to\infty}\frac{f(r) + 1}{g(r) + 1} = 1.
    \end{equation*}

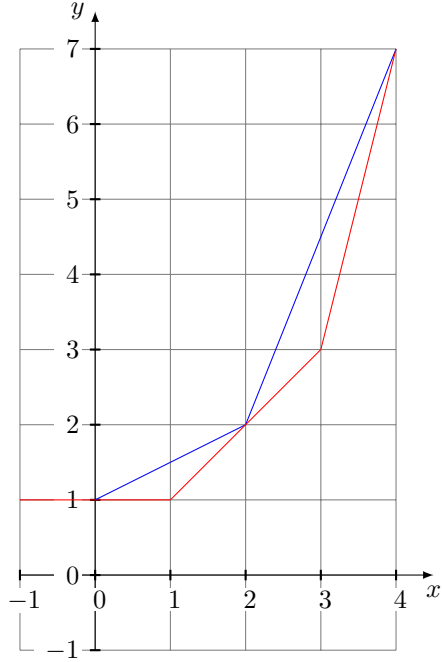
\begin{figure}[h!]
\centering 
  \begin{tikzpicture}
\def\a{1.5}
\tkzInit[xmax=4,ymax=7,xmin=-1,ymin=-1]
   \begin{scope}
        \tkzGrid
    \end{scope}
   \tkzAxeXY

\draw [domain=-1:1, color = red] plot(\x,1) node[right]{};

\draw [domain=0:2, color = blue] plot(\x,{(\a-1)*\x+1}) node[right]{};
\draw [domain=2:4, color = blue] plot(\x,{-8*\a*\a+12*\a-3+(4*\a*\a-5*\a+1)*\x}) node[right]{};

\draw [domain=1:3, color = red] plot(\x,{3-2*\a+(2*\a-2)*\x}) node[right]{};
\draw [domain=3:4, color = red] plot(\x,{-24*\a*\a+40*\a-15+(8*\a*\a-12*\a+4)*\x}) node[right]{};

\end{tikzpicture}    
  \caption{Functions $f$ and $g$ when $\psi(P,h) =\frac{2}{3}$ and $\Psi(P,h)=1$.}
  \label{fig:limsup}
\end{figure}

    In a similar way one can obtain any values for $\psi(P,h)$ and $\Psi(P,h)$ on the interval $[0,1]$ as long as $\psi(P,h)\leq \Psi(P,h)$. See Figure \ref{fig:limsup} for the case when $\psi(P,h) =\frac{2}{3}$ and $\Psi(P,h)=1$. 

The next corollary follows directly from Theorem \ref{theorem: second main theorem hypersurface new}.
\begin{apcorollary}\label{corollary: hyperplane second main theorem sum}
Let $q$ and $n$ be positive integers, $f : \R \to \mathbb{TP}^n$ a tropical holomorphic curve, $P_0,\ldots,P_q$ homogeneous tropical polynomials with tropical meromorphic coefficients and degrees $d_0,\ldots,d_q$, respectively.  
If for all the coefficients $a_j$ of all the polynomials $P_0,\ldots,P_q$ we have $N(r,a_j) = o(T_f(r))$, then
\begin{equation*}
    \sum_{j=0}^q \psi(P_j,f)T_f(r)\leq \sum_{j=0}^{q}\frac{1}{d_j}N\left(r,\frac{1_\circ}{P_j\circ f}\oslash\right) + o(T_f(r))\leq \sum_{j=0}^q \Psi(P_j,f)T_f(r).
\end{equation*}
\end{apcorollary}

If all of the coefficients of the tropical homogeneous polynomials have integer exponents and their coefficients are constants, then from Theorem \ref{theorem: hypersurface second main theorem old} it follows that
\begin{equation}\label{equation: hypersurface second main theorem old}
        \begin{aligned}
            &(q-M-\lambda)T_f(r)\\
            \leq &\sum_{j=M+1}^q \frac{1}{d_j} N\left(r,\ofrac{1_\circ}{P_j\circ f}\right) + o(T_f(r))\\
            \leq &(q-M)T_f(r),
        \end{aligned}
\end{equation}
where $r$ approaches infinity outside an exceptional set of zero upper density measure. Let $Q=q-M$. Then by renaming $P_{M+1},\ldots,P_q$ to $P_0,\ldots,P_{Q-1}$ we can write \eqref{equation: hypersurface second main theorem old} in the form
\begin{equation*}
    (Q-\lambda)T_f(r)
    \leq \sum_{j=0}^{Q-1} \frac{1}{d_j}N\left(r,\ofrac{1_\circ}{P_j\circ f}\right) + o(T_f(r))
    \leq QT_f(r),
\end{equation*}
where $r$ approaches infinity outside an exceptional set of zero upper density measure. On the other hand, Corollary \ref{corollary: hyperplane second main theorem sum} implies that
\begin{equation*}
    \sum_{j=0}^{Q-1} \psi(P_j,f)T_f(r)\leq \sum_{j=0}^{Q-1}\frac{1}{d_j}N\left(r,\frac{1_\circ}{P_j\circ f}\oslash\right) + o(T_f(r))\leq \sum_{j=0}^{Q-1} \Psi(P_j,f)T_f(r).
\end{equation*}
If $P_j\circ f$ is complete, then we can write it in the form
\begin{equation*}
    P_j(f)(x) = \max_{i_0+\cdots+i_n = d_j}\{c_{i_0,\ldots,i_n}+i_0f_0(x)+\cdots+i_nf_n(x)\},
\end{equation*}
where all the coefficients $c_j$ are real numbers. From this we can see that
\begin{equation*}
\begin{aligned}
    &\min_{i_0+\cdots+i_n = d_j}\{c_{i_0,\ldots,i_n}\} + d_j\max\{f_0(x),\ldots,f_n(x)\}\\
    &\leq P_j(f)(x)\\
    &\leq  \max_{i_0+\cdots+i_n = d_j}\{c_{i_0,\ldots,i_n}\} + d_j\max\{f_0(x),\ldots,f_n(x)\}.
\end{aligned}
\end{equation*}
 Thus when $P_j\circ f$ is complete, we have
\begin{equation*}\begin{aligned}
    \lim_{r\to\infty}\frac{\frac{1}{2}\left(P_j(f)(r)+P_j(f)(-r)\right)}{d_jT_f(r)} 
    = \lim_{r\to\infty}\frac{d_jT_f(r) + O(1)}{d_jT_f(r)} = 1.
\end{aligned}\end{equation*}
Therefore $Q-\lambda\leq \sum_{j=0}^{Q-1} \psi(P_j,f)T_f(r)$ and Lemma \ref{lemma: liminf limsup bounds} implies that $\sum_{j=0}^{Q-1} \Psi(P_j,f)\leq Q$.  In addition, any time we have $0 < \psi(P_j,f) < 1$ or $0 < \Psi(P_j,f) < 1$ for any $j$, we have $Q-\lambda< \sum_{j=0}^{Q-1} \psi(P_j,f)$ or $\sum_{j=0}^{Q-1} \Psi(P_j,f) < Q$ respectively. This shows that the bounds of Corollary \ref{corollary: hyperplane second main theorem sum} are tighter than in Theorem \ref{theorem: hypersurface second main theorem old}. By introducing a growth condition we can also obtain the following version of the second main theorem for tropical hypersurfaces.  

\begin{apcorollary} \label{corollary: second main thm with growth condition}
    Let $q$, $n$ and $d$ be positive integers such that $q > M$, where $M=\binom{n+d}{d}-1$. Let the tropical holomorphic curve $f:\R\to\mathbb{TP}^n$ be tropical algebraically nondegenerate. Assume that tropical hypersurfaces $V_{P_0},\ldots,V_{P_q}$ are defined by homogeneous tropical polynomials $P_0,\ldots,P_q$ with integer exponents and degrees $d_0,\ldots,d_q$, respectively, such that the least common multiple of $d_0,\ldots,d_q$ is equal to $d$. 
If $\lambda = \operatorname{ddg}(\{P_{M+1}\circ f,\ldots,P_{q}\circ f\})$ and
    \begin{equation*}
        \limsup_{r\to\infty}\frac{\log T_f(r)}{r} = 0,
    \end{equation*}
    then
    \begin{equation}\label{equation: corollary cao zheng with improved bounds statement}
        \begin{aligned}
            &\sum_{j=M+1}^q \psi(P_j,f)T_f(r)\\
            \leq &\sum_{j=0}^q\frac{1}{d_j}N\left(r,\ofrac{1_\circ}{P_j\circ f}\right) \\
            &- \frac{1}{d}N\left(r,\ofrac{1_\circ}{C_\circ(P_0^{\otimes\frac{d}{d_0}}\circ f,\ldots,P_{M}^{\otimes\frac{d}{d_{M}}}\circ f)}\right) + o(T_f(r))\\
            = &\sum_{j=M+1}^q \frac{1}{d_j}N\left(r,\ofrac{1_\circ}{P_j\circ f}\right) + o(T_f(r))\\
            \leq &\sum_{j=M+1}^q \Psi(P_j,f)T_f(r),
        \end{aligned}
    \end{equation}
    where $r$ approaches infinity outside an exceptional set of zero upper density measure.
\end{apcorollary}
\begin{proof}
    Theorem \ref{theorem: hypersurface second main theorem old} implies that
    \begin{equation*}
    \begin{aligned}
        &\sum_{j=0}^q\frac{1}{d_j}N\left(r,\ofrac{1_\circ}{P_j\circ f}\right) - \frac{1}{d}N\left(r,\ofrac{1_\circ}{C_\circ(P_0^{\otimes\frac{d}{d_0}}\circ f,\ldots,P_{M}^{\otimes\frac{d}{d_{M}}}\circ f)}\right) + o(T_f(r))\\
        &= \sum_{j=M+1}^q \frac{1}{d_j}N\left(r,\ofrac{1_\circ}{P_j\circ f}\right) + o(T_f(r))
    \end{aligned}
    \end{equation*}
    and by Corollary \ref{corollary: hyperplane second main theorem sum} we obtain
    \begin{equation*}
    \begin{aligned}
        \sum_{j=M+1}^q \psi(P_j,f)T_f(r)
        &\leq \sum_{j=M+1}^q \frac{1}{d_j}N\left(r,\ofrac{1_\circ}{P_j\circ f}\right) + o(T_f(r)) \\
        &\leq \sum_{j=M+1}^q \Psi(P_j,f)T_f(r).
    \end{aligned}
    \end{equation*}
\end{proof}
The growth condition from Theorem \ref{theorem: hypersurface second main theorem old} and Corollary \ref{corollary: second main thm with growth condition} cannot be dropped due to the equality in \eqref{equation: corollary cao zheng with improved bounds statement}. This can be seen in the following way. Let $f(x) = (0,e_2(x))$ and $P_0(x_0,x_1) = x_0\oplus x_1 = P_1(x_0,x_1) = P_2(x_0,x_1)$. Here $e_2(x)$ is the tropical hyper-exponential function. Now we have $n=1$, $d_0=d_1=d_2=d=1$, $M=1$ and $q = 2 > M$. We can see that
    \begin{equation*}
        \begin{aligned}
            &C_\circ(P_0^{\otimes\frac{d}{d_0}}\circ f,P_1^{\otimes\frac{d}{d_1}}\circ f)\\
            &=C_\circ(P_0\circ f,P_1\circ f)\\
            &=e_2(x)^++e_2(x+1)^+\\
            &=e_2(x)+2 e_2(x)\\
            &=3e_2(x).
        \end{aligned}
    \end{equation*}
    Now, on one hand,
    \begin{equation*}
    \begin{aligned}
        &\sum_{j=0}^2 N\left(r,\ofrac{1_\circ}{P_j\circ f}\right) - N\left(r,\ofrac{1_\circ}{C_\circ(P_0\circ f,P_1\circ f)}\right)\\
        & =3N\left(r,\ofrac{1_\circ}{e_2(x)}\right) - 3N\left(r,\ofrac{1_\circ}{e_2(x)}\right) \\
        & \equiv 0,\\
    \end{aligned}
    \end{equation*}
    but on the other hand,
    \begin{equation*}
        \sum_{j=M+1}^qN\left(r,\ofrac{1_\circ}{P_j\circ f}\right) =N\left(r,\ofrac{1_\circ}{e_2(x)}\right) = T\left(r,\ofrac{1_\circ}{e_2(x)}\right).
    \end{equation*}
    Since it is evidently true that $T\left(r,\ofrac{1_\circ}{e_2(x)}\right)\neq o(T_f(r))$, this shows that we cannot drop the growth condition. In this case also $\sum_{j=M+1}^q \psi(P_j,f) = \sum_{j=M+1}^q \Psi(P_j,f)= 1$, so changing the middle equality in \eqref{equation: corollary cao zheng with improved bounds statement} to an inequality will not help either.

In the case when the dimension of the tropical projective space is equal to one we have the following version of the second main theorem. 
\begin{apcorollary}\label{corollary: one dimensional holomorphic second main thm}
Let $f = [f_0:f_1] : \R \to \mathbb{TP}^1$ be a tropical holomorphic curve, $a = [a_0:a_1]$ a real constant and 
\begin{equation}\label{equation: dimension one polynomial}
    P(f)(x)=(a_0\otimes f_0(x))\oplus (a_1\otimes f_1(x)).
\end{equation}
Then
\begin{equation*}
\begin{aligned}
T(r,f) &= N\left(r,\frac{1_\circ}{P\circ f}\oslash\right) + O(1)
\end{aligned}
\end{equation*}
for all $a\in\R$.
\end{apcorollary}
\begin{proof}
    By Jensen formula we have
    \begin{equation*}
        \begin{aligned}
    &N\left(r,\frac{1_\circ}{P\circ f}\oslash\right) \\
    &= N\left(r,\frac{1_\circ}{P\circ f}\oslash\right) - N(r,P\circ f)\\
    &= \frac{1}{2}\sum_{\sigma=\pm 1}P(f)(\sigma r)-P(f)(0).
    \end{aligned}
    \end{equation*}
    Furthermore, 
    \begin{equation*}
        \min\{a_0,a_1\} + \|f(x)\|\leq (a_0\otimes f_0(x))\oplus (a_1\otimes f_1(x)) \leq \max\{a_0,a_1\} + \|f(x)\|,
    \end{equation*}
    which means that $P(f)(x) = \|f(x)\| + O(1)$. By combining the above equations we obtain
    \begin{equation*}
        T(r,f) = T_f(r) + O(1) = N\left(r,\frac{1_\circ}{P\circ f}\oslash\right) + O(1).
    \end{equation*}
\end{proof}
We can see a connection between $N\left(r,\ofrac{1_\circ}{P\circ f}\right)$ and $N\left(r,\ofrac{1_\circ}{f\oplus a}\right)$ in the following way. If $P(f)(x)=(a_0\otimes f_0(x))\oplus (a_1\otimes f_1(x))$ then by the tropical Jensen formula
\begin{equation*}
    \begin{aligned}
        &N\left(r,\ofrac{1_\circ}{P\circ f}\right)\\
        =& \frac{1}{2}\left(P(f)(r)+P(f)(-r)\right) + O(1)\\
        =& \frac{1}{2}\left((a_0\otimes f_0(r))\oplus (a_1\otimes f_1(r)) + (a_0\otimes f_0(-r))\oplus (a_1\otimes f_1(-r))\right) + O(1)\\
        =& \frac{1}{2}\left((a_0\oslash a_1 )\oplus (f_1(r)\oslash f_0(r)) + (a_0\oslash a_1)\oplus (f_1(-r)\oslash f_0(-r))\right)\\ &+\frac{1}{2}\left(f_0(r)+f_0(-r)\right) + O(1)\\
        =& \frac{1}{2}\left((a_1\oslash a_0 )\oplus (f_1(r)\oslash f_0(r))\right)\\ &+\frac{1}{2}\left( (a_1\oslash a_0)\oplus (f_1(-r)\oslash f_0(-r))+f_0(r)+f_0(-r)\right) + O(1)\\
        =& \frac{1}{2}\left(a\oplus f(r) + a\oplus f(-r) + f_0(r)+f_0(-r)\right) + O(1)\\
        =& N\left(r,\ofrac{1_\circ}{f\oplus a}\right)-N(r,f\oplus a) + N\left(r,\ofrac{1_\circ}{f_0}\right)-N(r,f_0) + O(1)\\
        =& N\left(r,\ofrac{1_\circ}{f\oplus a}\right)+ N\left(r,f\right)-N(r,f\oplus a)  + O(1).
    \end{aligned}
\end{equation*}
If $a$ is a tropical meromorphic function such that $m(r,a) = o(T(r,f))$ then similarly we can see that
\begin{equation*}
    N\left(r,\ofrac{1_\circ}{P\circ f}\right) = N\left(r,\ofrac{1_\circ}{f\oplus a}\right)+ N\left(r,f\right)-N(r,f\oplus a) + o(T(r,f)).
\end{equation*}
This means that Theorem \ref{theorem: second main theorem hypersurface new} implies Theorem \ref{theorem: second main theorem} for tropical meromorphic functions. 

\section{Tropical Nevanlinna inverse problem}
The defect for a tropical meromorphic function $f$ and a constant $a\in\R$,
\begin{equation*}
    \delta(a,f) := 1-\limsup_{r\to\infty}\frac{N\left(r,\ofrac{1_\circ}{f\oplus a}\right)}{T(r,f)},
\end{equation*}
was first considered by Laine and Tohge \cite{lainet:11}. Here we will consider the defect for tropical meromorphic functions $f$ and $a$. With Theorem \ref{theorem: second main theorem} we can see that if $m(r,a) = o(T(r,f))$ and $m\left(r,\ofrac{1_\circ}{f\oplus a}\right) = o(T(r,f))$, then
\begin{equation}\label{equation: defect different form}
\begin{aligned}
    \delta(a,f) 
    & = 1-\limsup_{r\to\infty}\frac{N\left(r,\ofrac{1_\circ}{f\oplus a}\right)}{T(r,f)}\\
    & = 1-\limsup_{r\to\infty}\frac{T(r,f)-N(r,f)+N(r,f\oplus a)}{T(r,f)}\\
    & = 1-\limsup_{r\to\infty}\left(1+\frac{N(r,f\oplus a)-N(r,f)}{T(r,f)}\right)\\
    & = -\limsup_{r\to\infty}\frac{N(r,f\oplus a)-N(r,f)}{T(r,f)}\\
    & = \liminf_{r\to\infty}\frac{N(r,f)-N(r,f\oplus a)}{T(r,f)}.\\
\end{aligned}
\end{equation}
From this formula we can see that in some sense the defect $\delta(a,f)$ measures how many of the poles of $f$ have a value greater than or equal to $a$. It is good to note that the term $N(r,f)-N(r,f\oplus a)$ also appears in Theorem \ref{theorem: second main theorem}. As we will see later on, this form of the defect will often be very useful.

In classical Nevanlinna theory the defect has the properties $0\leq\delta(a,f)\leq 1$ and 
\begin{equation*}
    \sum_{a\in\C\cup\{\infty\}}\delta(a,f)\leq 2.
\end{equation*}
In the tropical context the defect does have the property $0\leq\delta(a,f)\leq 1$, but it is possible, that  
\begin{equation*}
    \sum_{a\in\T}\delta(a,f) = \infty.
\end{equation*}
In fact
\begin{equation*}
    \sum_{a\in\T}\delta(a,f) < \infty
\end{equation*}
if and only if $\delta(a,f) = 0$ for all $a\in\T$. That is because $\delta(a,f)\leq \delta(b,f)$ for all $a\leq b$. In the tropical context it makes sense to consider a singe target value in the inverse problem since, unlike in the classical Nevanlinna theory, in tropical Nevanlinna theory the second main theorem can be stated for a single target value. Because of this we obtain the following versions of the inverse problem. 
\begin{aptheorem}\label{theorem: weaker inverse problem}
    For all $\delta\in[0,1]$ there exists a tropical rational function $f$ and a constant $a\in\R$ such that $\delta(a,f) = \delta$.
\end{aptheorem}

\begin{aptheorem}\label{theorem: stronger inverse problem}
    There exists a tropical meromorphic function $f$ such that for all $\delta\in[0,1]$ there exists $a\in\R$ such that $\delta(a,f) = \delta$.
\end{aptheorem}

Example \ref{example: weaker inverse solution} and Example \ref{example: stronger inverse solution} will prove the above theorems. In Theorem \ref{theorem: weaker inverse problem} it is interesting to note, that unlike in classical Nevanlinna theory, here we can solve the inverse problem with just a rational function. Theorem \ref{theorem: stronger inverse problem} cannot be solved with a rational function, but in some sense it resembles the classical inverse problem more, since we may choose multiple different targets for a single function.

\begin{example}\label{example: weaker inverse solution}
    Let $\alpha,\beta > 0$ and define
    \begin{equation*}
        f(x) = \max\left\{-\alpha\left|x-\frac{1}{\alpha}\right|+1,-\beta\left|x+\frac{2}{\beta}\right|+2\right\}.
    \end{equation*}
    Now $f$ has a pole of multiplicity $\alpha$ at $\frac{1}{\alpha}$ and a pole of multiplicity $\beta$ at $-\frac{2}{\beta}$. It is also true that $f(\frac{1}{\alpha}) = 1$ and  $f(-\frac{2}{\beta}) = 2$. Clearly now $m(r,f) = O(1)$ and thus by \eqref{equation: defect different form} we have
    \begin{equation*}
        \delta(x,f) = \begin{cases}
            1, &x\geq 2,\\
            \frac{\beta}{\alpha+\beta}, &1\leq x<2,\\
            0, &x<1.
        \end{cases}
    \end{equation*}
\end{example}
\begin{figure}
    \centering
    \begin{tikzpicture}
\def\a{0.4}
\def\b{1}
\tkzInit[xmax=5.7,ymax=2,xmin=-5,ymin=-1]
   \begin{scope}
    \end{scope}
   \tkzLabelX[orig=false,step=1,node font=\scriptsize]
   \tkzLabelY[orig=false,step=1,node font=\scriptsize]
   \tkzDrawX
   \tkzDrawY

\draw [domain=0:1/\a] plot(\x,\a*\x) node[right]{};
\draw [domain=1/\a:2/\a+1] plot(\x,-\a*\x+2) node[right]{};

\draw [domain=-2/\b:0] plot(\x,-\b*\x) node[right]{};
\draw [domain=-4/\b-1:-2/\b] plot(\x,\b*\x+4) node[right]{};








\end{tikzpicture}    
    \caption{$f(x)$ with $\alpha = \frac{2}{5}$ and $\beta=1$.}
    \label{fig:inverse_problem_rational}
\end{figure}
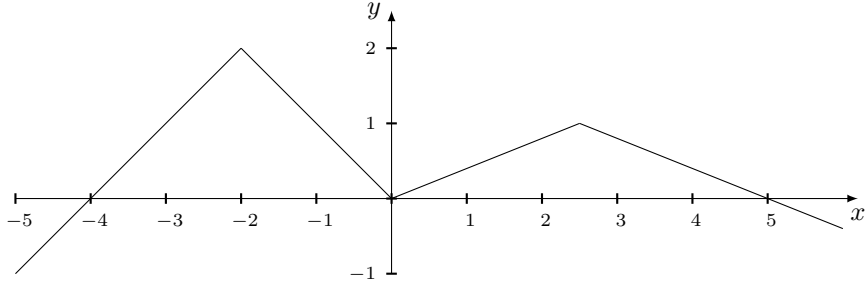

In general if we think of $\delta(x,f):\R\to[0,1]$ as a function of $x$ for some tropical meromorphic function $f$ we can say that $\delta(x,f)$ is a non-decreasing function. That is because in the form \eqref{equation: defect different form} we can see that with bigger and bigger values of $a$, more and more poles are flattened, and thus $N(r,f\oplus a)$ becomes smaller and smaller. The target value $a$ cannot create poles, so therefore $N(r,f\oplus a)$ cannot increase by making $a$ bigger. We also have the following result concerning the defect.
\begin{aptheorem}\label{theorem: defect piecewise constant}
Let $f$ be a tropical meromorphic function. 
Define the function $\delta(x,f):\R\to[0,1]$ as follows
\begin{equation*}
\begin{aligned}
    \delta(x,f) 
    & = 1-\limsup_{r\to\infty}\frac{N\left(r,\ofrac{1_\circ}{f\oplus x}\right)}{T(r,f)}.
\end{aligned}
\end{equation*}
If $x_0\in\R$ is not an accumulation point of $\{f(a):\omega_f(a)<0\}$, then $\delta(x,f)$ is locally constant at $x_0$, i.e. there exists $a<b$ such that $x_0\in[a,b]$ and $\delta(x,f)$ is constant for all $x\in[a,b]$.
\end{aptheorem}
\begin{proof}
If $f$ does not have any poles, then $\delta(x,f)\equiv 0$. If all of the poles of $f$ have the same value $c$ then for all $x_0 < c$ we have $\delta(x,f) = 0$ for all $x\leq x_0$ and for all $x_0\geq c$ we have $\delta(x,f) = \delta(c,f)$ for all $x\geq x_0$. So in that case $\delta(x,f)$ is locally constant for all $x\in\R$. Therefore we can assume that $f$ has at least two poles with different values. 
 
 If $x_0<\inf\{f(a):\omega_f(a)<0\}$, then $\delta(x,f) = 0$ for all $x\leq x_0$. Similarly if $x_0\geq s:=\sup\{f(a):\omega_f(a)<0\}$, then $\delta(x,f) = \delta(s,f)$ for all $x\geq x_0$. If $\inf\{f(a):\omega_f(a)<0\}\leq x_0<\sup\{f(a):\omega_f(a)<0\}$ and $x_0$ is not an accumulation point of $\{f(a):\omega_f(a)<0\}$, then we can find two poles with values $c_1,c_2\in\R$ such that $x_0\in[c_1,c_2]$ and there is no pole of $f$ with value $c$ such that $c\in[c_1,c_2]$. Now by \eqref{equation: defect different form} we can see that on the interval $[c_1,c_2]$ the function $\delta(x,f)$ is constant. 
\end{proof}

Sometimes it is possible that for a tropical meromorphic function $f$ the defect $\delta(x,f)$ is locally constant at the accumulation points of $\{f(a):\omega_f(a)<0\}$. For example, if $\inf\{f(a):\omega_f(a)<0\}$ or $\sup\{f(a):\omega_f(a)<0\}$ is an accumulation point, then by the proof of Theorem \ref{theorem: defect piecewise constant} we can see that $\delta(x,f)$ is still locally constant at those points. Later on in Example \ref{example: stronger inverse solution} we will see that defect is not always locally constant at some accumulation points.

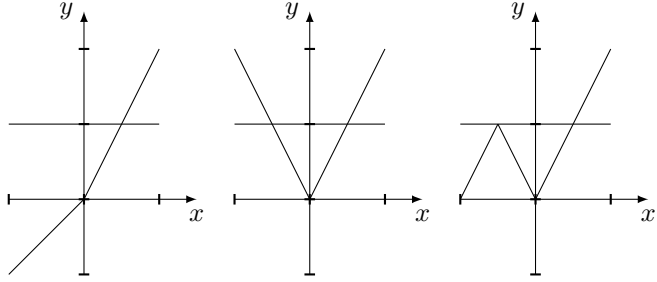
\begin{figure}[h!]
    \centering
    \begin{tabular}{c c c}
        \begin{tikzpicture}

\tkzInit[xmax=1,ymax=2,xmin=-1,ymin=-1]
   \begin{scope}
    \end{scope}
   \tkzDrawX
   \tkzDrawY

\draw [domain=(0:1] plot(\x,2*\x) node[right]{};
\draw [domain=-1:0] plot(\x,\x) node[right]{};
\draw [domain=-1:1] plot(\x,1) node[right]{};



\end{tikzpicture}     &  \begin{tikzpicture}

\tkzInit[xmax=1,ymax=2,xmin=-1,ymin=-1]
   \begin{scope}
    \end{scope}
   \tkzDrawX
   \tkzDrawY

\draw [domain=(0:1] plot(\x,2*\x) node[right]{};
\draw [domain=-1:0] plot(\x,-2*\x) node[right]{};
\draw [domain=-1:1] plot(\x,1) node[right]{};



\end{tikzpicture}     & \begin{tikzpicture}

\tkzInit[xmax=1,ymax=2,xmin=-1,ymin=-1]
   \begin{scope}
    \end{scope}
   \tkzDrawX
   \tkzDrawY

\draw [domain=0:1] plot(\x,2*\x) node[right]{};
\draw [domain=-0.5:0] plot(\x,-2*\x) node[right]{};
\draw [domain=-1:-0.5] plot(\x,2*\x+2) node[right]{};
\draw [domain=-1:1] plot(\x,1) node[right]{};



\end{tikzpicture}    \\
    \end{tabular}
    \caption{Different interactions between a root and the target value.}
    \label{fig:root flattening behaviour}
\end{figure}

We can see that when thought of as a function of $x$, the defect $\delta(x,f)$ seems to change value only at the values of the poles of $f$ even though $N\left(r,\ofrac{1_\circ}{f\oplus x}\right)$ is counting roots of $f$ and not poles. The interaction that the target value has with poles is more straightforward than with roots. When the target value goes above a value of a pole, it is flattened completely.  With roots, depending on the situation, when the target value goes above it we get one or two new roots such that the total multiplicity of the roots stays the same. Only when the target value goes over a pole, some contribution from the roots is lost. In fact the lost contribution is exactly equal to the multiplicity of the pole, which was flattened. This is evident when the defect is written in the form \eqref{equation: defect different form}, since the target value can only create roots but not poles. In Figure \ref{fig:root flattening behaviour} you can see from left to right the situations, when one root is created, two roots are created and when we lose contribution of the root equal to the multiplicity of the pole.

For a tropical meromorphic function $f$, if $\{f(a):\omega_f(a)<0\}$ does not have any accumulation points, then $\delta(x,f)$ is a piecewise constant function by Theorem \ref{theorem: defect piecewise constant}. In the next example we will construct a tropical meromorphic function $f$ such that the defect $\delta(x,f)$ is piecewise constant function of $x$ with infinitely many distinct values. 
\begin{example}
    Define the sequence $(a_n)_{n\in\N}$ as follows. First make every other element equal to 1
    \begin{equation*}
        (a_n) = (1,\_,1,\_,1,\_,1,\_,1,\_,1,\_,1,\_,1,\_,1,\_,1,\_,1,\_,1,\_,1,\_,1,\_,\ldots).
    \end{equation*}
    Then in the remaining blank spaces, starting from 2, make every other number equal to 2
    \begin{equation*}
        (a_n) = (1,2,1,\_,1,2,1,\_,1,2,1,\_,1,2,1,\_,1,2,1,\_,1,2,1,\_,1,2,1,\_,\ldots).
    \end{equation*}
    Then every other blank space becomes 3 and so on
    \begin{align*}
        (a_n) &= (1,2,1,3,1,2,1,\_,1,2,1,3,1,2,1,\_,1,2,1,3,1,2,1,\_,1,2,1,3,\ldots),\\
        (a_n) &= (1,2,1,3,1,2,1,4,1,2,1,3,1,2,1,\_,1,2,1,3,1,2,1,4,1,2,1,3,\ldots),\\
        (a_n) &= (1,2,1,3,1,2,1,4,1,2,1,3,1,2,1,5,1,2,1,3,1,2,1,4,1,2,1,3,\ldots).
    \end{align*}
    In the first step of the construction of the sequence we put $a_{2n-1} = 1$ for all $n\in\N$. Next we put $a_{2(2n-1)} = 2 = a_{2n-1}+1$, then $a_{2^2(2n-1)} = 3 = a_{2(2n-1)} + 1$, $a_{2^3(2n-1)} = 4 = a_{2^2(2n-1)} + 1$  and in general in the $(k+1)$th step we put $a_{2^k(2n-1)} = k + 1 = a_{2^{k-1}(2n-1)} + 1$ for all $n\in\N$. From this we can see that $a_{2n}  =  a_n + 1$ for all $n\in\N$.  
    
    Then define $f(x) = \max_{n\in\N}\left\{-\left|x-b(n)\right|+a_n\right\}$, where
    \begin{equation*}
        b(n) = 2\sum_{k=1}^n a_k-a_n = \sum_{k=1}^na_k+\sum_{k=1}^{n-1}a_k.
    \end{equation*}
    
    Now $f(x)$ is a tropical meromorphic function such that it has a pole of multiplicity 2 and value $a_n$ at $b(n)$ for every $n\in\N$. In addition $f(b(n)\pm a_n) =0$ for every $n\in\N$.  

    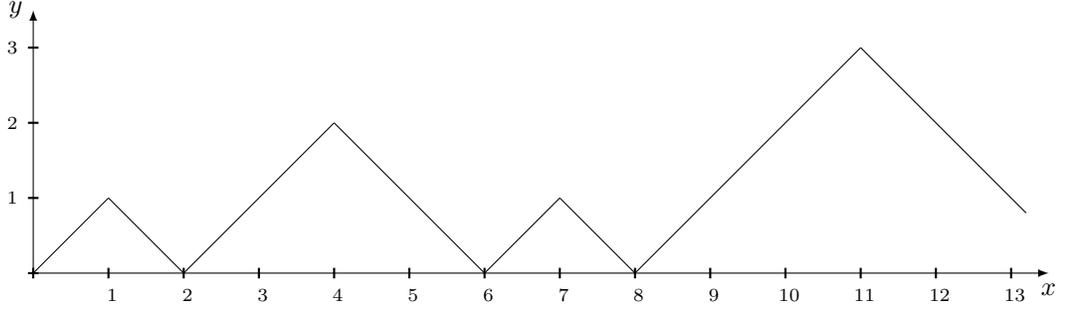
\begin{figure}[h]
\centering
  \begin{tikzpicture}
\def\xmin{0}
\tkzInit[xmax=13,ymax=3,xmin=\xmin,ymin=0]
   \begin{scope}
    \end{scope}
   \tkzLabelX[orig=false,step=1,node font=\scriptsize]
   \tkzLabelY[orig=false,step=1,node font=\scriptsize]
   \tkzDrawX
   \tkzDrawY

\draw [domain=\xmin:1] plot(\x,\x-2*1+1+1) node[right]{};
\draw [domain=1:2] plot(\x,-\x+2*1-1+1) node[right]{};

\draw [domain=2:4] plot(\x,\x-2*3+2*2) node[right]{};
\draw [domain=4:6] plot(\x,-\x+2*3) node[right]{};

\draw [domain=6:7] plot(\x,\x-2*4+2*1) node[right]{};
\draw [domain=7:8] plot(\x,-\x+2*4) node[right]{};

\draw [domain=8:11] plot(\x,\x-2*7+2*3) node[right]{};
\draw [domain=11:13.2] plot(\x,-\x+2*7) node[right]{};





\end{tikzpicture}    
  \caption{Function $f(x) = \max_{n\in\N}\{-|x-b(n)|+a_n\}$.}
  \label{fig:inverse}
\end{figure}
    
    Next we will examine the quantity $b(n)$. We can see that
    \begin{equation*}
        \begin{aligned}
            \sum_{k=1}^n a_k 
            &= \sum_{k=1}^{\left[\frac{n+1}{2}\right]}a_{2k-1} +\sum_{k=1}^{\left[\frac{n}{2}\right]}a_{2k} \\
            &= \left[\frac{n+1}{2}\right] +\sum_{k=1}^{\left[\frac{n}{2}\right]}(1+a_{k})\\
            &= \left[\frac{n+1}{2}\right] + \left[\frac{n}{2}\right] +\sum_{k=1}^{\left[\frac{n}{2}\right]}a_{k}\\
            &= n +\sum_{k=1}^{\left[\frac{n}{2}\right]}a_{k}\\
            &= n + \left[\frac{n}{2}\right] + \sum_{k=1}^{\left[\frac{n}{2^2}\right]}a_{k}\\
            &\vdots\\
            & =\sum_{k = 1}^\infty \left[\frac{n}{2^{k-1}}\right].
        \end{aligned}
    \end{equation*}
    The above sum is a finite sum, because for large enough $k$ we have $\left[\frac{n}{2^{k-1}}\right] = 0$. In fact we can see that $\left[\frac{n}{2^{k-1}}\right]\geq 1$ whenever $k\leq \log_2 n + 1$.
    Therefore 
    \begin{equation*}
        \sum_{k=1}^n a_k =\sum_{k = 1}^\infty \left[\frac{n}{2^{k-1}}\right]= \sum_{k = 1}^{[\log_2 n] + 1} \left[\frac{n}{2^{k-1}}\right].
    \end{equation*}
    We can see that
    \begin{equation*}\begin{aligned}
        \sum_{k = 1}^{[\log_2 n] + 1} \left[\frac{n}{2^{k-1}}\right] 
        &\leq \sum_{k = 1}^{[\log_2 n] + 1} \frac{n}{2^{k-1}} \\
        &= \left(2-\frac{1}{2^{[\log_2 n] + 1}}\right)n\\
        &\leq \left(2-\frac{1}{2^{\log_2 n + 1}}\right)n \\
        &= \left(2-\frac{1}{2n}\right)n = 2n-\frac{1}{2}.
    \end{aligned}\end{equation*}
    Similarly
    \begin{equation*}\begin{aligned}
        \sum_{k = 1}^{[\log_2 n] + 1} \left[\frac{n}{2^{k-1}}\right]
        &\geq \sum_{k = 1}^{[\log_2 n] + 1} \left(\frac{n}{2^{k-1}}-1\right)\\
        &=\left(2-\frac{1}{2^{[\log_2 n] + 1}}\right)n - [\log_2 n] - 1\\
        &\geq 2n - 1 - [\log_2 n] - 1 = 2n - [\log_2 n] - 2.
    \end{aligned}\end{equation*}
    So overall
    \begin{equation*}
        2n - [\log_2 n] - 2 \leq \sum_{k=1}^n a_k\leq 2n-\frac{1}{2}.
    \end{equation*}
    From this we obtain
    \begin{equation*}
        4n-[\log_2 n]-[\log_2 (n-1)]-6\leq b(n) \leq 4n - 3.
    \end{equation*}
    Define
    \begin{equation*}
        m(r) := \max\left\{n\in\N: b(n)\leq r\right\}.
    \end{equation*}
    We can see that
    \begin{equation*}
        r \leq b(m(r)+1) \leq 4m(r) + 1 \iff m(r) \geq \frac{r-1}{4}  
    \end{equation*}
    and
    \begin{equation*}\begin{aligned}
        r\geq b(m(r)) 
        \geq 4m(r) - [\log_2 m(r)]-[\log_2 (m(r)-1)]-6 
        \geq 2m(r) - 6, \\
    \end{aligned}\end{equation*}
    which means that $\frac{r-1}{4}\leq m(r)\leq \frac{r+6}{2}$ and thus $m(r) = O(r)$. We can also see that
    \begin{equation*}
        \sum_{n=1}^kb(n) \leq  \sum_{n=1}^k(4n-3) = 2k(k+1) - 3k = 2k^2-k
    \end{equation*}
    and
    \begin{equation*}\begin{aligned}
        \sum_{n=1}^kb(n) 
        &\geq \sum_{n=1}^k (4n-[\log_2 n]-[\log_2 (n-1)]-6) \\
        &= 2k(k+1) - 6k -\sum_{n=1}^k\left([\log_2 n]+[\log_2 (n-1)]\right) \\
        &= 2k^2 - 4k - \sum_{n=1}^k\left([\log_2 n]+[\log_2 (n-1)]\right).
    \end{aligned}\end{equation*}
    Therefore 
    \begin{equation*}
        \sum_{n=1}^kb(n) = 2k^2 + O(k\log k).
    \end{equation*}
    Then the counting function for $f$ is
    \begin{equation*}\begin{aligned}
        N(r,f) 
        &= \sum_{n=1}^{m(r)}\left(r-b(n)\right) \\
        &= rm(r) - \sum_{n=1}^{m(r)}b(n) \\
        &= rm(r) - 2m(r)^2 + O(m(r)\log m(r)) \\
        &= m(r)(r - 2m(r)) + O(r\log r).
    \end{aligned}\end{equation*}
    Also
    \begin{equation*}
        N(r,f)-N(r,f\oplus 1) = \sum_{n=1}^{\left[\frac{m(r)+1}{2}\right]}(r-b(2n-1))
    \end{equation*}
    and in general since $N(r,f)-N(r,f\oplus a)$ is counting only the poles of $f$ that have a value of $a$ or less, we can see that
    \begin{equation*}
        N(r,f)-N(r,f\oplus a) = \sum_{j=1}^k\sum_{n=1}^{\left[\frac{m(r)}{2^j}+\frac{1}{2}\right]}(r-b(2^{j-1}(2n-1))),
    \end{equation*}
    whenever $k$ is a positive integer and $a\in[k,k+1)$. Now
    \begin{equation*}\begin{aligned}
        &\sum_{n=1}^{\left[\frac{m(r)}{2^j}+\frac{1}{2}\right]}b(2^{j-1}(2n-1)) \\
        &= \sum_{n=1}^{\left[\frac{m(r)}{2^j}+\frac{1}{2}\right]} \left(4\cdot 2^{j-1}(2n-1) + O(\log n)\right)\\
        &= 4\cdot 2^{j-1}\left[\frac{m(r)}{2^j}+\frac{1}{2}\right]\left(\left[\frac{m(r)}{2^j}+\frac{1}{2}\right]+1\right) + O(r\log r)
    \end{aligned}\end{equation*}
    and
    \begin{equation*}\begin{aligned}
        &\sum_{j=1}^k\sum_{n=1}^{\left[\frac{m(r)}{2^j}+\frac{1}{2}\right]}b(2^{j-1}(2n-1))\\
        &\leq \sum_{j=1}^k4\cdot 2^{j-1}\left(\frac{m(r)}{2^j}+\frac{1}{2}\right)\left(\frac{m(r)}{2^j}+\frac{3}{2}\right) +  O(r\log r) \\
        &= 2m(r)^2\sum_{j=1}^k\frac{1}{2^j} + O(r\log r)\\
        &= 2m(r)^2\left(1-\frac{1}{2^k}\right) + O(r\log r)
    \end{aligned}\end{equation*}
    and similarly
    \begin{equation*}\begin{aligned}
        \sum_{j=1}^k\sum_{n=1}^{\left[\frac{m(r)}{2^j}+\frac{1}{2}\right]}b(2^{j-1}(2n-1))
        \geq 2m(r)^2\left(1-\frac{1}{2^k}\right) + O(r\log r).
    \end{aligned}\end{equation*}
    In a similar way we obtain
    \begin{equation*}
        \sum_{j=1}^k\sum_{n=1}^{\left[\frac{m(r)}{2^j}+\frac{1}{2}\right]}r = rm(r)\left(1-\frac{1}{2^k}\right).
    \end{equation*}
    By combining above equalities, we can see that
    \begin{equation*}\begin{aligned}
        \delta(a,f) 
        &= \limsup_{r\to\infty}\frac{N(r,f)-N(r,f\oplus a)}{N(r,f)} \\
        &= \limsup_{r\to\infty}\frac{\left(1-\frac{1}{2^k}\right)(r-2m(r))m(r)+ O(r\log r)}{(r-2m(r))m(r) + O(r\log r)} \\
        &= 1-\frac{1}{2^k}
    \end{aligned}\end{equation*}
    whenever $k$ is a positive integer and $a\in [k,k+1)$. This means that the defect is a piecewise constant function that attains a different value on each interval $[k,k+1)$ for $k\in\N$. 
\end{example}

\begin{aptheorem}\label{theorem: defect in terms of unintegrated counting function}
    Let $f:\R\to\R$  be a tropical meromorphic function and $a:\R\to\R$ a tropical entire function such that $m(r,a) = o(T(r,f))$ and $m\left(r,\ofrac{1_\circ}{f\oplus a}\right) = o(T(r,f))$. If $f$ has at least one pole, then
    \begin{equation*}
         \delta(a,f) \leq \limsup_{r\to\infty}\frac{n(r,f)-n(r,f\oplus a)}{n(r,f)}.
    \end{equation*}
    If in addition $m(r,f) = o(T(r,f))$, then
    \begin{equation*} \label{equation: unintegrated defect equality}
         \liminf_{r\to\infty}\frac{n(r,f)-n(r,f\oplus a)}{n(r,f)}\leq \delta(a,f) \leq \limsup_{r\to\infty}\frac{n(r,f)-n(r,f\oplus a)}{n(r,f)}.
    \end{equation*}
\end{aptheorem}
\begin{proof}
    If $N(r,f)-N(r,f\oplus a)\equiv 0$ then also $n(r,f)-n(r,f\oplus a) \equiv 0 $ and thus $\delta(a,f) = 0$. Assume that $N(r,f)-N(r,f\oplus a)\not\equiv 0$. Then we can see that $N(r,f)-N(r,f\oplus a)\to\infty$ and $N(r,f)\to\infty$ as $r\to\infty$ and by applying Theorem \ref{theorem: second main theorem} and the L'Hôpital's rule, we can see that
    \begin{equation*}\begin{aligned}
        \delta(a,f) &= \liminf_{r\to\infty}\frac{N(r,f)-N(r,f\oplus a)}{T(r,f)}\\
        &\leq \liminf_{r\to\infty}\frac{N(r,f)-N(r,f\oplus a)}{N(r,f)}\\
        &\leq \limsup_{r\to\infty}\frac{n(r,f)-n(r,f\oplus a)}{n(r,f)}.
    \end{aligned}\end{equation*}
    If $m(r,f) = o(T(r,f))$, then by the L'Hôpital's rule we obtain
    \begin{equation*}
        \begin{aligned}
        \delta(a,f) &= \liminf_{r\to\infty}\frac{N(r,f)-N(r,f\oplus a)}{T(r,f)}\\
        &= \liminf_{r\to\infty}\frac{N(r,f)-N(r,f\oplus a)}{N(r,f)}\\
        &\geq \liminf_{r\to\infty}\frac{n(r,f)-n(r,f\oplus a)}{n(r,f)}.
    \end{aligned}
    \end{equation*}
\end{proof}

The next example will show that sometimes at the accumulation points of $\{f(a):\omega_f(a)<0\}$ the defect is not locally constant. It also proves Theorem \ref{theorem: stronger inverse problem}, which is the inverse problem.
\begin{example}\label{example: stronger inverse solution}
    Let $F:[0,1]\to[0,1]$ be a non-decreasing function such that $F(0) = 0$ and $F(1) = 1$. By \cite[Chapter 2, Theorem 4.3]{kuipers:74} there exists a sequence $(a_n)\subset[0,1]$ such that
    \begin{equation*}
        \lim_{k\to\infty}\frac{\operatorname{card}(\{a_1,\ldots,a_k\}\cap[0,x])}{k} = F(x),
    \end{equation*}
    for all $x\in[0,1]$. Then define
    \begin{equation*}
        f(x) = \max_{n\in\N}\left\{-\left|x-b(n)\right|+a_n\right\},
    \end{equation*}
    where 
    \begin{equation*}
        b(n) = 2\sum_{k=1}^n a_k-a_n = \sum_{k=1}^n a_k+\sum_{k=1}^{n-1}a_k.
    \end{equation*}
    Define $m(r) := \max\left\{n\in\N: b(n)\leq r\right\}$ and
    \begin{equation*}
        I(a,x) = \begin{cases}
            1, &a\leq x,\\
            0, &a > x.
        \end{cases}
    \end{equation*}
    Now we can see that
    \begin{equation*}
        \begin{aligned}
            &\lim_{r\to\infty}\frac{n(r,f)-n(r,f\oplus x)}{n(r,f)}\\
            &= \lim_{r\to\infty}\frac{\sum_{n=1}^{m(r)}I(a_n,x)}{\sum_{n=1}^{m(r)}1}\\
            &= \lim_{r\to\infty}\frac{\operatorname{card}(\{a_1,\ldots,a_{m(r)}\}\cap[0,x])}{m(r)}
            = F(x),
        \end{aligned}
    \end{equation*}
    when $x\in[0,1]$. Theorem \ref{theorem: defect in terms of unintegrated counting function} now implies that
    \begin{equation*}
        \delta(x,f) = \begin{cases}
            1, & x\geq 1,\\
            F(x), & 0< x< 1,\\
            0, &x\leq 0.
        \end{cases}
    \end{equation*} 
\end{example}
If $F$ is chosen to be any increasing continuous function then $f$ satisfies Theorem~\ref{theorem: stronger inverse problem}.


Given a tropical hypersurface $V_P$ with a tropical homogeneous polynomial $P$ of degree $d$, Cao and Zheng \cite{caozheng:18} defined the defect as follows
\begin{equation*}
    \delta(V_P,f) := 1-\limsup_{r\to\infty}\frac{N\left(r,\ofrac{1_\circ}{P\circ f}\right)}{dT_f(r)}.
\end{equation*}

We aim to formulate the inverse problem for tropical hypersurfaces. To that end we will first give the following lemma. 
\begin{aplemma}\label{lemma: for any t there exists curve and polynomial}
    For all $t\in[0,1]$ there exists a tropical holomorphic curve $f:\R\to\mathbb{TP}^n$ and a homogeneous tropical polynomial $P$ of degree $d$ such that $\psi(P,f) = \Psi(P,f) = t$.
\end{aplemma}
\begin{proof}
    Define $f(x) = [t|x|:|x+1|:0:\cdots:0]$, where $t\in[0,1)$ and $P(x_0,\ldots,x_n) = x_0^{\otimes d}$. Now we have $P(f)(x) = dt|x|$ and
    \begin{equation*}
        T_f(r) = \frac{1}{2}(\max\{tr,r+1\}+\max\{tr,|r-1|\}) = r,
    \end{equation*}
    when $r\geq\frac{1}{1-t}$. Therefore
    \begin{equation*}
        \lim_{r\to\infty}\frac{\frac{1}{2}(P(f)(r)+P(f)(-r))}{dT_f(r)} = t = \psi(P,f) = \Psi(P,f),
    \end{equation*}
    for any $t\in[0,1)$. To get $\psi(P,f) = \Psi(P,f) = 1$ we can just choose 
    $$
    P(x_0,x_1,\ldots,x_n) = x_1^{\otimes d}
    $$
    instead. 
\end{proof}

We are now ready to formulate the inverse problem with tropical hypersurfaces. 
\begin{aptheorem}
    For all $\delta\in[0,1]$ there exists a tropical holomorphic curve $f:\R\to\mathbb{TP}^n$ and a homogeneous tropical polynomial $P$ of degree $d$ such that $\delta(V_P,f) = \delta$.
\end{aptheorem}
\begin{proof}
    By Lemma \ref{lemma: for any t there exists curve and polynomial} we can find a tropical holomorphic curve $f$ and a homogeneous tropical polynomial $P$ such that $\psi(P,f) = \Psi(P,f) = 1-\delta$. Theorem \ref{theorem: second main theorem hypersurface new} implies now that $\delta(V_P,f) = \delta$. 
\end{proof}

Cao and Zheng proved the following result regarding the defect relation, which can be seen as the tropical version of the Shiffman conjecture \cite{shiffman:77}.
\begin{aptheorem}[\cite{caozheng:18}]\label{theorem: cao zheng defect theorem}
    Let $q$, $n$ and $d$ be positive integers such that $q > M$, where $M=\binom{n+d}{d}-1$. Let the tropical holomorphic curve $f:\R\to\mathbb{TP}^n$ be tropical algebraically nondegenerate. Assume that tropical hypersurfaces $V_{P_0},\ldots,V_{P_q}$ are defined by homogeneous tropical polynomials $P_0,\ldots,P_q$ with degrees $d_0,\ldots,d_q$, respectively, such that the least common multiple of $d_0,\ldots,d_q$ is $d$. If $\lambda = \operatorname{ddg}(\{P_{M+1}\circ f,\ldots,P_{q}\circ f\})$ and
    \begin{equation*}
        \limsup_{r\to\infty}\frac{\log T_f(r)}{r} = 0,
    \end{equation*}
    then
    \begin{equation*}
        \sum_{j=0}^q\delta(V_{P_j},f)\leq M + 1 + \lambda, \text{ and } \sum_{j=M+1}^q\delta(V_{P_j},f)\leq \lambda.
    \end{equation*}
\end{aptheorem}

The following defect relation follows directly from Theorem \ref{theorem: second main theorem hypersurface new}.

\begin{aptheorem}\label{theorem: defect relation new}
Let $f : \R \to \mathbb{TP}^n$ be a non-constant tropical holomorphic curve and let $P_j:\mathbb{TP}^n\to\R(\not\equiv 0_\circ)$ be a tropical homogeneous polynomial of degree $d$ with tropical meromorphic coefficients. If for all coefficients $a_j$ of $P\circ f$ we have $N(r,a_j) = o(T_f(r))$, then 
\begin{equation}\label{equation: defect double ineq}
    1-\Psi(P,f)\leq\delta(V_{P},f)\leq 1-\psi(P,f).
\end{equation}
\end{aptheorem}
We can see that Theorem \ref{theorem: defect relation new} implies Theorem \ref{theorem: cao zheng defect theorem}. 
By Theorem \ref{theorem: defect relation new} we obtain
\begin{equation*}
    \sum_{j=0}^q\delta(V_{P_j},f)\leq q + 1 - \sum_{j = 0}^q \psi(P_j,f)
\end{equation*}
and
\begin{equation*}
    \sum_{j=M+1}^q\delta(V_{P_j},f)\leq q - M - \sum_{j = M+1}^q \psi(P_j,f).
\end{equation*}
We have already shown that if $P_j\circ f$ is complete, then $\psi(P_j,f) = 1$. 
This means that $\sum_{j = M+1}^q \psi(P_j,f)\geq q - M -\lambda$. From this immediately follows that
\begin{equation*}
    q - M - \sum_{j = M+1}^q \psi(P_j,f) \leq \lambda
\end{equation*}
and
\begin{equation*}
    q + 1 - \sum_{j = 0}^q \psi(P_j,f)\leq q + 1 - \sum_{j = M+1}^q \psi(P_j,f) \leq  M + 1 + \lambda.
\end{equation*}

Based on the tropical version of the Shiffman conjecture, Cao and Zheng also proposed the tropical version of Griffiths conjecture \cite{griffiths:72}.
\begin{apconjecture}[\cite{caozheng:18}]\label{conjecture: defect conjecture}
Let $q$, $n$ and $d$ be positive integers such that $q > M$, where $M=\binom{n+d}{d}-1$. Let the tropical holomorphic curve $f:\R\to\mathbb{TP}^n$ be tropical algebraically nondegenerate. Assume that tropical hypersurfaces $V_{P_0},\ldots,V_{P_q}$ are defined by homogeneous tropical polynomials $P_0,\ldots,P_q$ with degrees $d_0,\ldots,d_q$, respectively, such that the least common multiple of $d_0,\ldots,d_q$ is $d$. If $\lambda = \operatorname{ddg}(\{P_{M+1}\circ f,\ldots,P_{q}\circ f\})$ and
    \begin{equation*}
        \limsup_{r\to\infty}\frac{\log T_f(r)}{r} = 0,
    \end{equation*}
    then
    \begin{equation*}
        \sum_{j=0}^q\delta(V_{P_j},f) \leq \frac{n + 1 + \lambda}{d}.
    \end{equation*}
\end{apconjecture}
When $d=1$ we have $M=n$ and in that case the conjecture is true by Theorem~\ref{theorem: cao zheng defect theorem}. However, we can see that this conjecture is not true when $d>1$. By Lemma~\ref{lemma: for any t there exists curve and polynomial} we can find a holomorphic curve $f$ and a homogeneous tropical polynomial $P$ such that
\begin{equation}\label{equation: conjecture counterexample a = 0}
    \Psi(P,f) = \psi(P,f) = 0.
\end{equation}
By choosing $P_j = P$ for all $j = 0,\ldots,q$ we obtain 
\begin{equation*}
     q + 1 = \sum_{j=0}^q\delta(V_{P_j},f)
\end{equation*}
and $\lambda = q-M$. 
In this case it can be seen that
\begin{equation}\label{equation: conjecture disp eq 1}
    \frac{n + 1 + \lambda}{d} = \frac{n + q + 1 - M}{d} < q + 1,
\end{equation}
whenever
\begin{equation*}
    q>\frac{n - M}{d-1}-1.
\end{equation*}
Since $n\leq M$ in general, we can see that \eqref{equation: conjecture disp eq 1} is actually true for all values of $q$. This contradicts Conjecture \ref{conjecture: defect conjecture}, because we would have 
\begin{equation*}
    \frac{n + 1 + \lambda}{d}  < q + 1= \sum_{j=0}^q\delta(V_{P_j},f)\leq \frac{n + 1 + \lambda}{d}.
\end{equation*}
We can also see that Conjecture \ref{conjecture: defect conjecture} is not true even if you replace $\frac{n + 1 + \lambda}{d}$ with $\frac{M + 1 + \lambda}{d}$.

It is good to note that in the simplest case, when $n=1,d=1$, $a = [a_0:a_1]$ and $P(x_0,x_1) = a_0\otimes x_0\oplus a_1\otimes x_1$ and the tropical holomorphic curve $f = [f_0:f_1]$ is thought of as a tropical meromorphic function $f=f_1\oslash f_0$, then as we have seen before 
\begin{equation*}
    N\left(r,\ofrac{1_\circ}{P\circ f}\right) = N\left(r,\ofrac{1_\circ}{f\oplus a}\right) + N(r,f) - N(r,f\oplus a) + O(1).
\end{equation*}
Thus in this case the defect $\delta(V_P,f)$ cannot be seen as a generalization of $\delta(a,f)$. In fact, in this simplest case we always have $\delta(V_P,f) = 0$ by Theorem \ref{theorem: second main theorem}.

\section{Appendix}
\subsection{Inverse problem in classical Nevanlinna theory}
In this section of the appendix we will state the inverse problem in the classical Nevanlinna theory. In order to do that we need to give some definitions, starting with the Nevanlinna functions. The proximity function for a meromorphic function $f$ in classical Nevanlinna theory is defined as 
\begin{equation*}
    m(r,f) := \frac{1}{2\pi}\int_0^{2\pi}\log^+\!|f(re^{i\theta})|\,d\theta,
\end{equation*}
where $\log^+\!x = \max\{\log x, 0 \}$. The Nevanlinna counting function is defined as
\begin{equation*}
    N(r,f) := \int_0^r\frac{n(t,f)-n(0,f)}{t}\,dt + n(0,f)\log r,
\end{equation*}
where $n(r,f)$ counts the poles of $f$ in $|z|<r$ according to their multiplicities. The Nevanlinna characteristic function is then defined as 
\begin{equation*}
    T(r,f) := m(r,f) + N(r,f).
\end{equation*}

For more details about classical Nevanlinna theory see for example~\cite{cherryy:01,goldbergo:08,hayman:64}.

In order to state the inverse problem, we also need to define the defect and index of multiplicity for meromorphic functions. The defect for a meromorphic function $f$ is defined as
\begin{equation*}
    \delta(a,f) =1-\limsup_{r\to\infty}\frac{N\left(r,\frac{1}{f-a}\right)}{T(r,f)}
\end{equation*}
for $a\in\C$ and
\begin{equation*}
    \delta(a,f) := 1-\limsup_{r\to\infty}\frac{N\left(r,f\right)}{T(r,f)}
\end{equation*}
for $a  = \infty$. The index of multiplicity is defined as
\begin{equation*}
    \theta(a,f) := \liminf_{r\to\infty}\frac{N\left(r,\frac{1}{f-a}\right)-\overline{N}\left(r,\frac{1}{f-a}\right)}{T(r,f)},
\end{equation*}
for $a\in\C$ and
\begin{equation*}
    \theta(a,f) :=\liminf_{r\to\infty}\frac{N\left(r,f\right)-\overline{N}\left(r,f\right)}{T(r,f)},
\end{equation*}
for $a = \infty$. Here $\overline{N}(r,f)$ is the truncated counting function, which counts the poles of $f$ without taking multiplicity into account. Using the second main theorem one can prove that
\begin{equation*}
    0\leq \delta(a,f) + \theta(a,f) \leq 1
\end{equation*}
for any meromorphic function $f$ and a constant $a\in\C\cup\{\infty\}$ and
\begin{equation*}
    \sum_{a\in\C\cup\{\infty\}}\delta(a,f)+\theta(a,f)\leq 2.
\end{equation*}
The inverse problem in classical Nevanlinna theory is then as follows. For $1\leq i < N\leq\infty$ let sequences $\{\delta_i\},\ \{\theta_i\}$ of non-negative numbers be assigned such that
    \begin{equation*}
        0<\delta_i+\theta_i\leq 1
    \end{equation*}
    for all $1\leq i < N$ and
    \begin{equation*}
        \sum_i(\delta_i+\theta_i)\leq 2.
    \end{equation*}
    Let $\{a_i\},\ 1\leq i< N$ be a sequence of distinct complex numbers. Does there exists a meromorphic function $f$ such that
    \begin{equation*}
        \delta(a_i,f) = \delta_i,\,\theta(a_i,f) = \theta_i
    \end{equation*}
    for all $1\leq i < N$ and
    \begin{equation*}
        \delta(a,f) = \theta(a,f) = 0
    \end{equation*}
for all $a\not\in\{a_i\}$? Drasin  \cite{drasin:77} has shown using quasiconformal mappings that such a meromorphic function always exists.

\subsection{Tropical Nevanlinna theory}
The tropical semiring is defined as $\T = \R\cup\{-\infty\}$ where addition and multiplication are defined as
\begin{equation*}
    a\oplus b = \max\{a,b\}
\end{equation*}
and
\begin{equation*}
    a\otimes b = a + b. 
\end{equation*}
Additive and multiplicative neutral elements are $0_\circ = -\infty$ and $1_\circ = 0$. Here $\T$ is a semiring, because not all elements have an additive inverse element. For example there is no $x\in\T$ such that $2\oplus x = 0_\circ$. For this reason subtraction is not defined on the tropical semiring. Tropical division is defined as $a\oslash b = a-b$ and exponentiation as $a^{\otimes\alpha} = \alpha a$ for $\alpha\in\R$. For more details about tropical geometry see for example \cite{maclagans:15}.

Next we will define tropical meromorphic functions.
\begin{definition}[\cite{halburds:09,lainet:11}]
    A continuous piecewise linear function $f:\R\to\R$ is said to be tropical meromorphic. 
\end{definition}

We say that $x$ is a pole of $f$ if 
\begin{equation*}
    \omega_f(x) := \lim_{\varepsilon\to0^+}(f'(x+\varepsilon)-f'(x-\varepsilon)) < 0
\end{equation*}
and a root if $\omega_f(x)>0$. The multiplicity of a root or a pole of a tropical meromorphic function $f$ at $x$ is $\tau_f(x) := |\omega_f(x)|$. A tropical meromorphic function is called a tropical entire function, if it has no poles. A tropical polynomial is a tropical entire function with finitely many roots and a tropical rational function is a tropical meromorphic function that has finitely many roots and poles. Tropical polynomials can be written in the form
\begin{equation*}
    \bigoplus_{n=1}^k c_n \otimes x^{\otimes t_n} = \max_{n=1}^k\{c_n + t_n x\}
\end{equation*}
and tropical entire functions can be written in the form
\begin{equation*}
    \bigoplus_{n=1}^\infty c_n \otimes x^{\otimes t_n} = \max_{n=1}^\infty\{c_n + t_n x\},
\end{equation*}
where $c_n\in\T$ and $t_n\in\R$ \cite{korhonenlainetohge:15}. Every tropical meromorphic function $h$ can be written in the form
\begin{equation}\label{equation: tropical meromorphic as a quotient of entire functions}
    h = \ofrac{f}{g},
\end{equation}
where $f$ and $g$ are tropical entire functions which do not share any roots \cite[Proposition 3.3]{korhonent:16AM}. If $f$ and $g$ in \eqref{equation: tropical meromorphic as a quotient of entire functions} are polynomials, then $h$ is a tropical rational function. 

Next we will define the tropical Nevanlinna functions. The tropical proximity function is defined as 
\begin{equation*}
    m(r,f) = \frac{f(r)^++f(-r)^+}{2},
\end{equation*}
where $f(x)^+ = \max\{f(x),0\}$. The tropical counting function is defined as 
\begin{equation*}
    N(r,f) = \frac{1}{2}\int_0^rn(t,f)dt = \frac{1}{2}\sum_{|b_\nu|<r}\tau_f(b_\nu)(r-|b_\nu|),
\end{equation*}
where $n(r,f)$ counts the poles of $f$ in $(-r,r)$ according to their multiplicities. The tropical Nevanlinna characteristic function is then defined as 
\begin{equation*}
    T(r,f) := m(r,f) + N(r,f).
\end{equation*}
The tropical Nevanlinna characteristic is a positive, convex, continuous, non-decreasing piecewise linear function of $r$ \cite[Lemma 3.2]{halburds:09}.  

The order and hyper-order of tropical meromorphic function $f$ are defined as
\begin{align*}
    \rho(f) &= \limsup_{r\to\infty}\frac{\log T(r,f)}{\log r}\\
    \rho_2(f) &= \limsup_{r\to\infty}\frac{\log\log  T(r,f)}{\log r}.
\end{align*}

As a special case of the tropical Poisson-Jensen formula (Theorem \ref{theorem: poisson jensen old}) we have the tropical Jensen formula when $x=0$
\begin{equation*}
\begin{aligned}
    f(0) &=\frac{1}{2}(f(r)+f(-r))
    - \frac{1}{2}\sum_{|a_\mu|<r}\tau_f(a_\mu)(r-|a_\mu|)+\frac{1}{2}\sum_{|b_\nu|<r}\tau_f(b_\nu)(r-|b_\nu|)\\
    &=m(r,f)-m\left(r,\ofrac{1_\circ}{f}\right) + N(r,f) -N\left(r,\ofrac{1_\circ}{f}\right)\\ 
    &= T(r,f) -T\left(r,\ofrac{1_\circ}{f}\right).
\end{aligned}
\end{equation*}

We can list some basic properties for the Nevanlinna functions \cite[Lemma 3.2]{korhonenlainetohge:15}.

\begin{enumerate}[(i)]
    \item If $f\leq g$, then $m(r,f)\leq m(r,g)$.
    \item Given a positive real number $\alpha$, then
    \begin{align*}
        &m(r,f^{\otimes\alpha}) = \alpha m(r,f),\\
        &N(r,f^{\otimes\alpha}) = \alpha N(r,f),\\
        &T(r,f^{\otimes\alpha}) = \alpha T(r,f).\\   
    \end{align*}
    \item Given tropical meromorphic functions $f,g$, then
    \begin{align*}
        &m(r,f\otimes g)\leq m(r,f)+m(r,g),\\
        &N(r,f\otimes g)\leq N(r,f)+N(r,g),\\
        &T(r,f\otimes g)\leq T(r,f)+T(r,g),\\
    \end{align*}
    \item and similarly,
    \begin{align*}
        &m(r,f\oplus g)\leq m(r,f)+m(r,g),\\
        &N(r,f\oplus g)\leq N(r,f)+N(r,g),\\
        &T(r,f\oplus g)\leq T(r,f)+T(r,g).
    \end{align*}
\end{enumerate}

There also exists a tropical counterpart to the lemma on the logarithmic derivative. Cao and Zheng proved the following version of the tropical lemma on the logarithmic derivative.
\begin{aptheorem}[\cite{caozheng:18}]
    Let $c\in\R\!\setminus\! \{0\}$. If $f$ is a tropical meromorphic function on $\R$ with
    \begin{equation*}
        \limsup_{r\to\infty}\frac{\log T(r,f)}{r} = 0,
    \end{equation*}
    then
    \begin{equation*}
        m\left(r,\ofrac{f(x+c)}{f(x)}\right) = o(T(r,f)),
    \end{equation*}
    where $r$ runs to infinity outside of a set of zero upper density measure $E$, i.e.
    \begin{equation*}
        \overline{\operatorname{dens}}E=\limsup_{r\to\infty}\frac{1}{r}\int_{E\cap[1,r]}dt = 0.
    \end{equation*}
\end{aptheorem}

Next we will introduce tropical hyper-exponential functions. Tropical hyper-exponential functions are reminiscent of hyper-exponential functions $\exp(z^c)$ over the usual algebra.
\begin{definition}[\cite{lainet:11}]
    Let $\alpha,\beta$ be real numbers such that $|\alpha|>1$ and $|\beta|<1$. Define functions $e_\alpha(x)$ and $e_\beta(x)$ on $\R$ by
    \begin{equation*}
        e_\alpha(x) := \alpha^{[x]}(x-[x])+\sum_{j=-\infty}^{[x]-1}\alpha^j=\alpha^{[x]}\left(x-[x]+\frac{1}{\alpha-1}\right)
    \end{equation*}
    and 
    \begin{equation*}
        e_\beta(x) := \sum_{j=[x]}^{\infty}\beta^j-\beta^{[x]}(x-[x])=\beta^{[x]}\left(\frac{1}{1-\beta}-x+[x]\right).
    \end{equation*}
\end{definition}
The tropical hyper-exponential functions have the following properties.
\begin{aplemma}[\cite{lainet:11}]
    Let $\alpha,\beta$ be real numbers such that $|\alpha|>1$ and $|\beta|<1$. The function $e_\alpha(x)$ is tropical meromorphic on $\R$ satisfying
    \begin{enumerate}[(i)]
        \item $e_\alpha(m) = \frac{\alpha^m}{\alpha-1}$ for each $m\in\Z$;
        \item $e_\alpha(x) = x+\frac{1}{\alpha-1}$ for any $x\in[0,1)$;
        \item the functional equation $y(x+1) = y(x)^{\otimes \alpha}$ on the whole $\R$.
    \end{enumerate}
    Similarly the function $e_\beta(x)$ is tropical meromorphic on $\R$ satisfying
    \begin{enumerate}[(i)]
        \item $e_\beta(m) = \frac{\beta^m}{1-\beta}$ for each $m\in\Z$;
        \item $e_\beta(x) = -x+\frac{1}{1-\beta}$ for any $x\in[0,1)$;
        \item the functional equation $y(x+1) = y(x)^{\otimes \beta}$ on the whole $\R$.
    \end{enumerate}
\end{aplemma}
There is a connection between $e_\alpha(x)$ and $e_\beta(x)$. Suppose $\alpha\neq \pm 1$. Then $e_\alpha(-x) = \frac{1}{\alpha}e_{\frac{1}{\alpha}}(x)$ for all $x\in\R$.

The tropical hyper-exponential function is of infinite order and hyper-order 1 \cite[Proposition 8.5]{lainet:11}.
\subsection{Tropical linear algebra}
First we will introduce tropical matrices. The operations of addition $\oplus$ and multiplication $\otimes$ for the $(n+1)\times (n+1)$ matrices $A=(a_{ij})$ and $B=(b_{ij})$ are defined as
\begin{equation*}
    A\oplus B = (a_{ij}\oplus b_{ij})
\end{equation*}
and
\begin{equation*}
    A\otimes B = \left(\bigoplus_{k=0}^na_{ik}\otimes b_{kj}\right),
\end{equation*}
respectively. The matrix $A$ is called regular if it contains at least one element different from $0_\circ$ in each row. The tropical determinant $|A|_\circ$ of $A$ is defined as \cite{korhonent:16AM,yoeli:61} 
\begin{equation*}
    |A|_\circ = \bigoplus a_{0\pi(0)}\otimes a_{1\pi(1)}\otimes\cdots\otimes a_{n\pi(n)},
\end{equation*}
where the tropical sum is taken over all permutations $\{\pi(0),\pi(1),\ldots,\pi(n)\}$ of $\{0,1,\ldots,n\}$. 

Now we can define the tropical Casorati determinant. 
Let $g(x)$ be a tropical entire function, $n\in\N$ and $c\in\R\setminus\{0\}$. For brevity we denote
\begin{equation*}
    g(x)\equiv g,\quad g(x+c)\equiv \overline g,\quad g(x+2c)\equiv \overline{\overline{g}}\quad \text{and} \quad g(x+nc)\equiv \overline g^{[n]}.
\end{equation*}
Now the tropical Casorati determinant of tropical entire functions $g_0,\ldots,g_n$ is defined by
\begin{equation*}
    C_\circ(g_0,g_1,\ldots,g_n)=\bigoplus \overline g_0^{[\pi(0)]}\otimes \overline g_1^{[\pi(1)]}\otimes\cdots\otimes \overline g_n^{[\pi(n)]},
\end{equation*}
where the tropical sum is taken over all permutations $\{\pi(0),\pi(1),\ldots,\pi(n)\}$ of $\{0,1,\ldots,n\}$. The tropical Casoratian has the following properties.
\begin{aplemma}[\cite{korhonent:16AM}]
    If $g_0,g_1,\ldots,g_n$ and $h$ are tropical entire functions, then
    \begin{enumerate}[(i)]
        \item $C_\circ(g_0,g_1,\ldots,g_i,\ldots,g_j,\ldots,g_n) = C_\circ(g_0,g_1,\ldots,g_j,\ldots,g_i,\ldots,g_n)$ for all $i,j\in\{0,\ldots,n\}$ such that $i\neq j$.
        \item $C_\circ(1_\circ,g_1,\ldots,g_n)\geq C_\circ(\overline g_1,\ldots,\overline g_n)$.
        \item $C_\circ(0_\circ,g_1,\ldots,g_n) = 0_\circ$.
        \item $C_\circ(g_0\otimes h,g_1\otimes h,\ldots,g_n\otimes h) = h\otimes\overline h\otimes\cdots\otimes\overline{h}^{[n]}\otimes C_\circ(g_0,g_1,\ldots,g_n)$.
    \end{enumerate}
\end{aplemma}

Next we will define tropical linear combinations and linear independence. 
\begin{definition}[\cite{korhonent:16AM}]
    If $g_0,\ldots,g_n$ are tropical meromorphic functions and $a_0,\ldots,a_n\in\T$, then
    \begin{equation*}
        f = \bigoplus_{\nu = 0}^n a_\nu\otimes g_\nu = \bigoplus_{i = 0}^j a_{k_i}\otimes g_{k_i}
    \end{equation*}
    is called a tropical linear combination of $g_0,\ldots,g_n$ over $\T$, where the index set $\{k_0,\ldots,k_j\}\subset\{0,\ldots,n\}$ is such that $a_{k_i}\in\R$ for all $i\in\{0,\ldots,j\}$, while $a_\nu=0_\circ$ if $\nu\not\in\{k_0,\ldots,k_j\}$.
\end{definition}
With the tropical linear combinations we can define linear independence in the way of Gondran and Minoux.
\begin{definition}[\cite{gondranm:79,gondranm:84,korhonent:16AM}]
    Tropical meromorphic functions $f_0,\ldots,f_n$ are linearly dependent (respectively independent) in the Gondran-Minoux sense if there exist (respectively do not exist) two disjoint subsets $I$ and $J$ of $K:=\{0,\ldots,n\}$ such that $I\cup J= K$ and
    \begin{equation*}
        \bigoplus_{i\in I}\alpha_i\otimes f_i = \bigoplus_{j\in J}\alpha_j\otimes f_j,
    \end{equation*}
    where the constants $\alpha_0,\ldots,\alpha_n\in\T$ are not all equal to $0_\circ$. 
\end{definition}
\begin{example}
    Let $g_0(x) = x-5, g_1(x) = |x|,g_2(x) = 1-x, g_3(x) = 2x$. Then $g_0,g_1,g_2,g_3$ are linearly dependent in the Gondran-Minoux sense because
    \begin{equation*}
        5\otimes g_0(x)\oplus (-1)\otimes g_2(x) 
        = 1_\circ\otimes g_1(x)\oplus 0_\circ\otimes g_3(x).
    \end{equation*}
    If you do not include the function $g_1$, then $g_0,g_2,g_3$ are linearly independent in the Gondran-Minoux sense.
\end{example}

Next we will define the notion of degeneracy.  
\begin{definition}[\cite{korhonent:16AM}]
    Let $G=\{g_0,\ldots,g_n\}(\neq \{0_\circ\})$ be a set of tropical entire functions linearly independent in the Gondran-Minoux sense, and let
    \begin{equation*}
        \mathcal{L}_G=\operatorname{span}\langle g_0,\ldots,g_n\rangle = \left\{\bigoplus_{k=0}^n a_k\otimes g_k : (a_0,\ldots,a_n)\in\T^{n+1}\right\}
    \end{equation*}
    to be their linear span. The collection $G$ is called the spanning basis of $\mathcal{L}_G$. The shortest length of the representation of $f\in\mathcal{L}_G\setminus\{0_\circ\}$ is defined by
    \begin{equation*}
        \ell(f)=\min\left\{j\in\{1,\ldots,n+1\}:f=\bigoplus_{i=1}^j a_{k_i}\otimes g_{k_i} \right\},
    \end{equation*}
    where $a_{k_i}\in\R$ with integers $0\leq k_1<k_2<\cdots<k_j\leq n$, and the dimension of $\mathcal{L}_G$ is
    \begin{equation*}
        \dim(\mathcal{L}_G) = \max\{\ell(f):f\in\mathcal{L}_G\setminus\{0_\circ\}\}.
    \end{equation*}
\end{definition}
\begin{definition}[\cite{korhonent:16AM}]
    Let $G=\{g_0,\ldots,g_n\}(\neq \{0_\circ\})$ be a set of tropical entire functions linearly independent in the Gondran-Minoux sense, and let $f$ be a tropical linear combination of $g_0,\ldots,g_n$. If $\ell(f) = n+1$, then $f$ is said to be complete. 
\end{definition}
\begin{definition}[\cite{korhonent:16AM}]
    Let $G=\{g_0,\ldots,g_n\}$ be a set of tropical entire functions, linearly independent in the Gondran-Minoux sense, and let $Q\subset\mathcal{L}_G$ be a collection of tropical linear combinations of $G$ over $\T$. The degree of degeneracy of $Q$ is defined to be 
    \begin{equation*}
        \operatorname{ddg}(Q) = \operatorname{card}(\{f\in Q:\ell(f)<n+1\}).
    \end{equation*}
    If $\operatorname{ddg}(Q)=0$, then $Q$ is called non-degenerate.
    
\end{definition}
\begin{example}
    Let $g_0(x) = x,g_1(x) = 2x, g_2(x) = 3x$. The tropical polynomials $g_0,g_1,g_2$ are linearly independent in the Gondran-Minoux sense so we can ask if their tropical linear combinations are complete or not. Consider the following tropical linear combinations 
    \begin{align*}
        f_0 &= 0_\circ\otimes g_0\oplus 5\otimes g_1\oplus 1_\circ\otimes g_2\\
        f_1 &= 1_\circ\otimes g_0\oplus 5\otimes g_1\oplus 1_\circ\otimes g_2\\
        f_2 &= 1_\circ\otimes g_0\oplus 1_\circ\otimes g_1\oplus 1_\circ\otimes g_2.
    \end{align*}
    Any tropical linear combination where at least one coefficient is $0_\circ$ is not complete, which means that $f_0$ is not complete. The tropical linear combination $f_2$ can be written in the form $f_2 = 1_\circ\otimes g_0\oplus 1_\circ\otimes g_2$, which means that it is not complete either. The tropical linear combination $f_1$ is of the form
    \begin{equation*}
        f_1(x) =\begin{cases}
            x,&x<-5,\\
            2x,&-5\leq x\leq 5,\\
            3x, &x>5.
        \end{cases}
    \end{equation*}
    Therefore $f_1$ is complete and $\operatorname{ddg}(\{f_0,f_1,f_2\}) = 2$.
\end{example}

\subsection{Tropical holomorphic curves and hypersurfaces}
We will now define the tropical projective space $\mathbb{TP}^n$. The space $\mathbb{TP}^n$ is defined in the following way. Let the equivalence relation $\sim$ be defined so that
\begin{equation*}
    (a_0,a_1,\ldots,a_n)\sim(b_0,b_1,\ldots,b_n)
\end{equation*}
if and only if
\begin{equation*}
    (a_0,a_1,\ldots,a_n) = \lambda\otimes (b_0,b_1,\ldots,b_n):= (\lambda\otimes b_0,\lambda\otimes b_1,\ldots,\lambda\otimes b_n)
\end{equation*}
for some $\lambda\in\R$. We denote by $[a_0:a_1:\cdots:a_n]$ the equivalence class of $(a_0,a_1,\ldots,a_n)$. The tropical projective space is now defined as the quotient space of $\T^{n+1}\setminus\{\mathbf{0}_\circ\}$ by the equivalence relation $\sim$, where $\mathbf 0_\circ = (0_\circ,\ldots,0_\circ)$ is the zero element of $\T^{n+1}$. The one dimensional tropical projective space $\mathbb{TP}^1$ can be identified with the completed tropical semiring $\T\cup\{\infty\}$ by the map
\begin{align*}
    &[1_\circ : a]\mapsto a\oslash 1_\circ = a, \quad a\in\T,\\
    &[0_\circ:a] \mapsto a\oslash 0_\circ = \infty,\quad a\in\R.
\end{align*}
We can now define the holomorphic curve.
\begin{definition}[\cite{korhonent:16AM}]
    Let $[a_0:\cdots:a_n]\in\mathbb{TP}^n$ be the equivalence class of $(a_0,\ldots,a_n)\in\T^{n+1}\setminus\{\mathbf{0}_\circ\}$, and let
    \begin{equation*}
        f = [g_0:\cdots:g_n]:\R\to\mathbb{TP}^n
    \end{equation*}
    be a tropical holomorphic map where $g_0,\ldots,g_n$ are tropical entire functions that do not have any roots which are common to all of them.
\end{definition}
Denote
\begin{equation*}
    \mathbf{f} = (g_0,\ldots,g_n):\R\to\R^{n+1}.
\end{equation*}
Then $\mathbf{f}$ is called a reduced representation of the tropical holomorphic curve $f$ in $\mathbb{TP}^n$. Next we will define the Cartan characteristic function for tropical holomorphic curves.
\begin{definition}[\cite{korhonent:16AM}]
    If $f:\R\to\mathbb{TP}^n$ is a tropical holomorphic curve with a reduced representation $\mathbf{f} = (g_0,\ldots,g_n)$, then
    \begin{equation*}
        T_{\mathbf{f}}(r) = \frac{1}{2}\left(\|\mathbf{f}(r)\|+\|\mathbf{f}(-r)\|\right), \quad \|\mathbf{f}(x)\| = \max\{g_0(x),\ldots,g_n(x)\}
    \end{equation*}
    is said to be the tropical Cartan characteristic function of $f$.
\end{definition}
It has been shown that the tropical Cartan characteristic does not depend on the reduced representation \cite[Proposition 4.3]{korhonent:16AM} so we will denote the tropical Cartan characteristic function $T_{\mathbf{f}}(r)$ simply as $T_f(r)$. 

Any tropical meromorphic function $f$ can always be represented as a quotient of two entire functions which do not share any common roots $f = h\oslash g$ \cite[Proposition 3.3]{korhonent:16AM}. If $f$ is now represented also in the form $f=[g,h]$, then in the case of one dimensional tropical projective space $\mathbb{TP}^1$ it has been shown that $T_f(r)$ is up to a constant equal to $T(r,f)$ \cite[Proposition 4.4]{korhonent:16AM}. The tropical Cartan characteristic then shares many of the properties of the tropical Nevanlinna characteristic. 

We will now define tropical hypersurfaces. 
\begin{definition}[\cite{caozheng:18}]
    Let $P$ be a homogeneous tropical polynomial in $n$-dimensional tropical projective space $\mathbb{TP}^n$. The set of roots of $P$ is called a tropical hypersurface. Denote the tropical hypersurface by $V_P$. 
\end{definition}
Set $M=\binom{n+d}{d}-1$. Now the composition of a tropical holomorphic curve $f=[f_0:f_1:\cdots:f_n]$ and a tropical homogeneous polynomial in $\mathbb{TP}^n$ of dimension $d$ can be written in the form
\begin{equation*}
    P\circ f = \bigoplus_{I_i\in \mathcal{J}_d}c_{I_i}\otimes f^{I_i} =\bigoplus_{i = 0}^M c_{I_i}\otimes f^{I_i} ,
\end{equation*}
where $\mathcal{J}_d$ is the set of all $I_i = (i_0,i_1,\ldots,i_n)\in\N_0^{n+1}$ such that $i_0+i_1+\cdots+i_n = d$ and $f^{I_i} := f_0^{\otimes i_0}\otimes\cdots\otimes f_n^{\otimes i_n}$. We can see that $P\circ f$ is a tropical algebraic combination. From this we have the following definition.
\begin{definition}[\cite{caozheng:18}]
    Tropical meromorphic functions $f_0,\ldots,f_n$ are algebraically dependent (respectively independent) in the Gondran-Minoux sense if $f^{I_0},\ldots,f^{I_M}$ are linearly dependent (respectively independent) in the Gondran-Minoux sense.
\end{definition}
\begin{definition}[\cite{caozheng:18}]
    Let $G=\{f_0,\ldots,f_n\}(\neq \{0_\circ\})$ be a set of tropical entire functions algebraically independent in the Gondran-Minoux sense, and let
    \begin{equation*}
        \hat{\mathcal{L}}_G=\operatorname{span}\langle f^{I_0},f^{I_1},\ldots,f^{I_M}\rangle = \left\{\bigoplus_{k=0}^n a_k\otimes f^{I_k} : (a_0,\ldots,a_n)\in\T^{M+1}\right\}
    \end{equation*}
    be their linear span. The collection $G$ is called the spanning basis of $\hat{\mathcal{L}}_G$. The shortest length of the representation of $f\in\hat{\mathcal{L}}_G\setminus\{0_\circ\}$ is defined as
    \begin{equation*}
        \hat{\ell}(F)=\min\left\{j\in\{1,\ldots,M+1\}:F=\bigoplus_{i=1}^j a_{k_i}\otimes f^{I_k} \right\},
    \end{equation*}
    where $a_{k_i}\in\R$ with integers $0\leq k_1<k_2<\cdots<k_j\leq n$, and the dimension of $\hat{\mathcal{L}}_G$ is
    \begin{equation*}
        \dim(\hat{\mathcal{L}}_G) = \max\{\hat{\ell}(F):F\in\hat{\mathcal{L}}_G\setminus\{0_\circ\}\}.
    \end{equation*}
\end{definition}
\begin{definition}[\cite{caozheng:18}]
    Let $G=\{f_0,\ldots,f_n\}(\neq \{0_\circ\})$ be a set of tropical entire functions linearly independent in the Gondran-Minoux sense, and let $F$ be a tropical algebraic combination of $f_0,\ldots,f_n$. If $\hat{\ell}(F) = M+1$, then $F$ is said to be complete. 
\end{definition}
Algebraic nondegeneracy is defined in an analogous way to classical algebraic geometry. 
\begin{definition}[\cite{caozheng:18}]
    Let $f=[f_0:f_1:\cdots:f_n]:\R\to\mathbb{TP}^n$ be a tropical holomorphic curve. If for any tropical hypersurface $V_P$ in $\mathbb{TP}^n$ defined by a homogeneous tropical polynomial $P$ in $\R^{n+1}$, $f(\R)$ is not a subset of $V_P$ then we say that $f$ is tropical algebraically nondegenerate. 
\end{definition}
There is a link between algebraic nondegeneracy and algebraic independence in the Gondran-Minoux sense as follows. 
\begin{aplemma}[\cite{caozheng:18}]
    A tropical holomorphic curve $f:\R\to\mathbb{TP}^n$ with reduced representation $f=(f_0,f_1,\ldots,f_n)$ is tropical algebraically nondegenerate if and only if $f_0,\ldots,f_n$ are algebraically independent in the Gondran-Minoux sense.
\end{aplemma}

\bibliographystyle{amsplain}
\bibliography{database}

\end{document}